\documentclass[12pt]{amsart}

\usepackage{amsmath,mathtools,amssymb,amsthm,amsfonts,amscd,mathrsfs}
\usepackage{dsfont,bigdelim,braket,physics,bm,latexsym}
\usepackage{graphicx,wrapfig,pgfplots}
\usepackage[dvipsnames]{xcolor}
\usepackage[normalem]{ulem}
\usepackage{comment,multirow,longtable,enumerate,enumitem}
\usepackage[all]{xy}
\usepackage{tikz,tikz-cd}
\usetikzlibrary{positioning,arrows.meta,fit,calc,backgrounds}
\pgfdeclarelayer{background}
\pgfdeclarelayer{foreground}
\pgfsetlayers{background,main,foreground}
\usepackage[most]{tcolorbox}
\tcbuselibrary{breakable,skins,theorems}
\tcbset{mybox/.style={colback=white,colframe=green!35!black,fonttitle=\bfseries,breakable=true}}
\usepackage[final]{showkeys}
\usepackage{url,xurl}
\usepackage{amsmath,amssymb,graphicx,tikz}
\usetikzlibrary{arrows.meta,positioning,fit,calc}
\usepackage[colorlinks,linkcolor=red,anchorcolor=blue,citecolor=blue]{hyperref}
\hypersetup{pdftitle={An Exposition of Kobayashi's Conjecture},pdfauthor={Masataka Iwai}}
\usepackage{letltxmacro}
\LetLtxMacro\manuscriptcite\cite
\RenewDocumentCommand{\cite}{o m}{%
  \ifmmode\text{\IfNoValueTF{#1}{\manuscriptcite{#2}}{\manuscriptcite[#1]{#2}}}%
  \else\IfNoValueTF{#1}{\manuscriptcite{#2}}{\manuscriptcite[#1]{#2}}\fi}
\setlist[enumerate]{label=$(\arabic*)$}
\newtheorem{thm}{Theorem}[section]

\newtheorem{cor}[thm]{Corollary}

\newtheorem{prop}[thm]{Proposition}

\newtheorem{conjecture}[thm]{Conjecture}

\newtheorem{lem}[thm]{Lemma}

\theoremstyle{definition} 
\newtheorem{defn}[thm]{Definition}

\newtheorem{question}[thm]{Question}

\theoremstyle{remark}
\newtheorem{rem}[thm]{Remark}

\newtheorem{setup}[thm]{Setup}

\newtheorem{claim}[thm]{Claim}

\numberwithin{equation}{section}

\DeclareMathOperator{\Ric}{Ric}

\newcommand{\underalign}[2]{%
  \quad
  \underset{\mathclap{\raisebox{-1.6ex}{\ensuremath{\scriptstyle #1}}}}{#2}
  \quad
}

\newcommand{\R}{\mathbb{R}}

\newcommand{\C}{\mathbb{C}}

\newcommand{\D}{\mathbb{D}}

\title{An Exposition of Kobayashi's canonical ampleness conjecture}
\author{Masataka Iwai}
\address{Department of Mathematics, Graduate School of Science, the University of Osaka,
1-1, Machikaneyama-cho, Toyonaka, Osaka 560-0043, Japan.}
\email{{\tt masataka@math.sci.osaka-u.ac.jp}}
\email{{\tt masataka.math@gmail.com}}
\date{\today, version 0.01}

\begin{document}

\begin{abstract}
We give an exposition of Kobayashi's canonical ampleness conjecture, that is, any compact K\"ahler manifold $X$ with no nonconstant holomorphic map $\C \to X$ has an ample canonical bundle. 
This result was announced in item 51 of the OpenAI mathematics catalogue~\cite{OpenAIMathCatalogue} (see \cite{OpenAI2026Canonical} for the proof).
We give a largely self-contained exposition. However, we replace some proofs of arguments that were difficult to follow or that the author did not understand.
We also include an approach using the Kähler--Ricci flow.
\end{abstract}
\maketitle
\setcounter{tocdepth}{3}
\tableofcontents

\section{Introduction}
\subsection{Kobayashi's canonical ampleness conjecture}

In this note, we give an exposition of the following formulation of Kobayashi's canonical ampleness conjecture.

\begin{thm}
\label{conj:intro-kobayashi}
Let $X$ be a compact K\"ahler manifold.
If $X$ is Brody hyperbolic, that is, if every holomorphic map
$\mathbb C\to X$ is constant, then the canonical bundle $K_X$ is ample.
\end{thm}

Notice that, for compact complex manifolds, Brody hyperbolicity and Kobayashi hyperbolicity are equivalent 
(cf.~\cite[Brody criterion~1.9]{Demailly2012Hyperbolic}).
Under the assumption that there exists a K\"ahler metric of negative holomorphic sectional curvature,
Wu--Yau proved the ampleness of $K_X$ for projective manifolds
\cite[Theorem~2]{Ch8WuYau}.
Tosatti--Yang extended this result to compact K\"ahler manifolds
and also proved nefness under nonpositive holomorphic sectional curvature
\cite[Theorem~1.1, Corollary~1.2]{IntroTosattiYang}.
Nomura subsequently gave alternative proofs of these results using the K\"ahler--Ricci flow
\cite[Theorems~1.1--1.2]{Nomura}.

In this paper, we present the argument in the OpenAI manuscript
\cite{OpenAI2026Canonical}, which assumes Brody hyperbolicity,
with an emphasis on the second variation of holomorphic discs.
We aim to make the proofs largely self-contained and include explanations even for elementary arguments.

\subsection{Outline of the proof}
Let $h_0$ be a fixed K\"ahler metric,
let $h$ be a K\"ahler metric, and let $q$ be a smooth closed real $(1,1)$-form.
We write $\omega_h$ for the K\"ahler form associated with a K\"ahler metric $h$.
For a holomorphic disc $f:\mathbb D\to X$,
the functionals in Definition~\ref{defn-APJ-funcional} are
\[
d(f):=|f'(0)|_{h_0}^2,\quad
P_q(f):=\frac2\pi\int_{\mathbb D}\log\frac1{|z|}
\,q(f'(z),f'(z))\,dS(z)
\]
\[
A_h(f):=\frac1\pi\int_{\mathbb D}h(f'(z),f'(z))\,dS(z)
\]
Here, $dS$ is the Euclidean area element on the disc.
If $X$ is Brody hyperbolic, then we can take a constant $D_0>0$, depending only on $X,h_0$, as in Lemma~\ref{lem:brody-bound},
so that every holomorphic disc $f$ satisfies
\[
d(f)=|f'(0)|_{h_0}^2\le D_0
\]
The following flowchart summarizes the proof. 

\newpage

\begin{center}
\resizebox{\linewidth}{!}{%
\begin{tikzpicture}[
  box/.style={
    rounded corners=3mm,
    draw=blue!55!black,
    fill=blue!5,
    text width=11.4cm,
    inner sep=4mm,
    align=center,
    font=\large,
    line width=.8pt
  },
  greenbox/.style={
    box,
    draw=green!45!black,
    fill=green!5
  },
  edge/.style={
    -{Stealth[length=3mm]},
    line width=1.1pt,
    draw=blue!60!black
  },
  note/.style={
    font=\large,
    fill=white,
    inner sep=1.5mm,
    align=center,
    text width=10cm
  }
]

% First row: align the top edges
\node[box,anchor=north west] (max) at (0,0)
{\textbf{(1) Maximizing discs and integral lower bounds}\\[2mm]
Find a maximizing disc $f : \D \to X$ with $J_\delta(f):=
d(f)
+P_q(f)
-\delta A_h(f)$
\\[2mm]
Estimate $
A_q(f)\ge-D_0, P_q(f)\ge-D_0
$\\[3mm]
Lemma~\ref{lem:brody-bound},
Proposition~\ref{prop:disc-maximizer},
Corollary~\ref{lem:signed-reward-source-bounds};
\cite[Lemmas~3.1, 3.3, Propositions~3.2, 3.4]
{OpenAI2026Canonical}
};

\node[box,anchor=north west] (variation)
at ([xshift=8mm]max.north east)
{\textbf{(2) Frames and holomorphic deformations}\\[2mm]
Find a good frame $(B_1, \ldots, B_n)$ of $f^{*}T_X$ with $\displaystyle
h(B_i,B_j)|_{\partial\mathbb D}=\delta_{ij}
$\\[2mm]
Compute 2nd variations of $d(f), P_q(f),  A_h(f)$ via deformations with initial velocity $zB_j$
\\[3mm]
Theorem~\ref{thm:disc-frame-variation-hessian},
Theorem~\ref{thm:disc-center-fixed-deformation},
Theorem~\ref{thm:disc-second-variation}; \\
cf.~\cite[\S4, Proposition~5.2]{OpenAI2026Canonical};
cf.~\cite[Lemma~3.3, Proposition~3.4]{DuGuoXie2026}
};

% Determine the bounding box of the first row
\node[
  fit=(max)(variation),
  inner sep=0pt,
  outer sep=0pt
] (toprow) {};

\node[font=\Large\bfseries,above=8mm of toprow]
{Flowchart of the proof};

% Second row: place below the first row
\node[
  box,
  text width=24.4cm,
  fill=blue!12,
  below=12mm of toprow
] (main)
{\textbf{(3) The basic inequality from local maximality and second variation}\\[3mm]
Using a maximizing disc $f : \D \to X $ with $x=f(0)$ and setting $q = \operatorname{Ric}(h)$, we obtain \\[2mm]
$\displaystyle
n\exp\!\left(
\frac{
P_{\operatorname{Ric}(h)}(f)
-\log\frac{\omega_h^n}{\omega_{h_0}^n}(x)
}{n}
\right)
+\frac1{2\pi}\int_0^{2\pi}
R_h(f(e^{\sqrt{-1}\theta}))\,d\theta
\le
\delta\bigl(n-A_{\operatorname{Ric}(h)}(f)\bigr)
$\\[2mm]
Theorem~\ref{thm:signed-disc-determinant-trace},
Proposition~\ref{prop:ricci-reward-master-inequality};
cf.~\cite[Corollary~5.3]{OpenAI2026Canonical}
};

% Third row: align the top edges
\node[greenbox,anchor=north west] (nef)
at ([yshift=-24mm]main.south west)
{\textbf{(4) Nefness via a disc envelope}\\[2mm]
$\displaystyle
\sup_{f\in\mathcal A_x}
P_{\operatorname{Ric}(h_0)}(f)\le C
\quad\Longrightarrow\quad K_X\text{ is nef}
$\\[2mm]
Use a disc envelope \\[3mm]
Corollary~\ref{cor:uniform-ricci-poisson-bound},
Theorem~\ref{thm:canonical-nef-envelope};\\
\cite[Proposition~6.1]{OpenAI2026Canonical}; 
cf.~\cite[Propositions~6.2, 6.4]{OpenAI2026Canonical}
};

\node[box,anchor=north east] (volume)
at ([yshift=-24mm]main.south east)
{\textbf{(5) A volume lower bound from scalar curvature $R_h$}\\[2mm]
$\displaystyle
R_h\ge-B
\quad\Longrightarrow\quad
\frac{\omega_h^n}{\omega_{h_0}^n}
\ge e^{-D_0}\left(\frac nB\right)^n
$\\[2mm]
Theorem~\ref{thm:scalar-volume-by-finite-prediction};
cf.~\cite[Proposition~7.1]{OpenAI2026Canonical}
};

% Determine the bounding box of the third row
\node[
  fit=(nef)(volume),
  inner sep=0pt,
  outer sep=0pt
] (bottomrow) {};

% Fourth row: place below the third row
\node[
  greenbox,
  text width=24.4cm,
  below=24mm of bottomrow
] (ma)
{\textbf{(6) Bigness via the Monge-Ampere equation}\\[3mm]
For $\varepsilon>0$, Take a metric 
$\displaystyle
\omega_{h_\varepsilon}
:=\varepsilon\omega_{h_0}-\operatorname{Ric}(h_0)
+\sqrt{-1}\,\partial\bar\partial u_\varepsilon>0,
\qquad
\omega_{h_\varepsilon}^n=e^{u_\varepsilon}\omega_{h_0}^n
$\\[3mm]
Apply the volume lower bound in (5), integrate, and let $\varepsilon \to 0$:\\[2mm]
$\displaystyle
(2\pi)^nc_1(K_X)^n
\ge e^{-D_0}\int_X\omega_{h_0}^n>0
$\\[2mm]
Together with nefness, this shows that $K_X$ is big\\[3mm]
Theorem~\ref{thm:canonical-big-volume}; \\
\cite[Proposition~8(i)]{Ch8WuYau};
cf.~\cite[Theorem~8.1, Proposition~8.2]{OpenAI2026Canonical};
\cite[Corollary~6.19]{DemaillyAnalyticMethods};
\cite[Lemma~9.1]{OpenAI2026Canonical}
};

% Arrows
\draw[edge]
(max.south) -- (main.north -| max.south);

\draw[edge]
(variation.south) -- (main.north -| variation.south);

\draw[edge]
(main.south -| nef.north) -- (nef.north)
node[midway,note]
{Compare with maximizing discs\\$\delta\downarrow0$};

\draw[edge]
(main.south -| volume.north) -- (volume.north)
node[midway,note]
{Use the lower bounds in (1) and let $\delta\downarrow0$};

\draw[edge,draw=green!45!black]
(nef.south) -- (ma.north -| nef.south)
node[midway,note,text=green!45!black]
{Construct metrics using nefness and Monge-Ampere Equation};

\draw[edge]
(volume.south) -- (ma.north -| volume.south)
node[midway,note,text=blue!60!black]
{Apply the volume lower bound to Monge-Ampere Equation};

\end{tikzpicture}%
}
\end{center}

We explain the proof following (1)--(6) in the flowchart.

\begin{enumerate}
\item In Proposition~\ref{prop:disc-maximizer}, for any $x\in X$ and $\delta >0$, we show that among holomorphic discs $f : \D \to X$ of finite area ($A_h(f)<+\infty$) with fixed center $f(0)=x$, the functional
\[
J_\delta(f):=
d(f)
+P_q(f)
-\delta A_h(f)
\]
attains its maximum. We also show that every maximizing disc is continuous on the closed disc and belongs to $W^{1,p}$ for $2<p<4$.
Furthermore, Corollary~\ref{lem:signed-reward-source-bounds} gives
%最大性を縮小円板$f_r(z):=f(rz)$および定数円板と比較すると、
\[
A_q(f)\ge-d(f)\ge-D_0,\qquad
P_q(f)\ge-D_0
\]
%を得る（系\ref{lem:signed-reward-source-bounds}）。
These lower bounds will be used in (3) and (5).

\item For any holomorphic disc $f : \D \to X$ that is continuous on the closed disc and belongs to $W^{1,p}$ for $2<p<4$, Theorem~\ref{thm:disc-frame-variation-hessian}(2) gives
a holomorphic frame $B=(B_1,\ldots,B_n)$ such that
\[
B_j\in W^{1,p}(\mathbb D,f^*TX),\qquad
h(B_i,B_j)|_{\partial\mathbb D}=\delta_{ij}
\]
By Theorem~\ref{thm:disc-center-fixed-deformation}, for each $j=1,\ldots, n$ there exists a deformation, holomorphic also in the parameter $s$,
\[
F_j:\{|s|<\varepsilon_j\}\longrightarrow W^{1,p}(\mathbb D,X),
\quad
F_j(0,z)=f(z),\quad F_j(s,0)=x,\quad
\left.\partial_sF_j(s,z)\right|_{s=0}=zB_j(z)
\]
Each $F_j(s,\cdot)$ is a holomorphic disc.
(This is the hardest part in this proof.)

In Theorem~\ref{thm:disc-second-variation}, we prove
\[
\begin{aligned}
\left.\partial_s\partial_{\bar s}d(F_j(s,\cdot))\right|_{s=0}
&=|B_j(0)|_{h_0}^2,\\
\left.\partial_s\partial_{\bar s}P_q(F_j(s,\cdot))\right|_{s=0}
&=\frac1{2\pi}\int_0^{2\pi}
q(B_j,B_j)(e^{\sqrt{-1}\theta})\,d\theta,\\
\left.\partial_s\partial_{\bar s}A_h(F_j(s,\cdot))\right|_{s=0}
&=\frac1\pi\int_{\mathbb D}
\partial_z\partial_{\bar z}
\bigl(|z|^2h(B_j,B_j)\bigr)\,dS.
\end{aligned}
\]

\item For a maximizing disc $f : \D \to X$ obtained in (1),
the local maximality of $J_\delta$ implies that its second variation is nonpositive. Substituting (2), we obtain
\[
\underbrace{\sum_{j=1}^n|B_j(0)|_{h_0}^2}
_{d\text{: sum of second variations}}
+
\underbrace{\frac1{2\pi}\int_0^{2\pi}
(\operatorname{tr}_h q)(f(e^{\sqrt{-1}\theta}))\,d\theta}
_{P_q\text{: sum of second variations}}
\le
\delta\underbrace{\bigl\{n-A_{\operatorname{Ric}(h)}(f)\bigr\}}
_{A_h\text{: sum of second variations}}
\]
We estimate the central term from below using  Green's formula.
Setting $q=\operatorname{Ric}(h)$ and writing $R_h = \operatorname{tr}_h\operatorname{Ric}(h)$ for the scalar curvature, we obtain
\[
\begin{aligned}
\underbrace{
n\exp\!\left(
\frac{P_{\operatorname{Ric}(h)}(f)
-\log\frac{\omega_h^n}{\omega_{h_0}^n}(x)}n
\right)}
_{d\text{ term}}
+
\underbrace{\frac1{2\pi}\int_0^{2\pi}
R_h(f(e^{\sqrt{-1}\theta}))\,d\theta}
_{P_{\operatorname{Ric}(h)}\text{ term}}
\le
\delta\underbrace{\bigl\{n-A_{\operatorname{Ric}(h)}(f)\bigr\}}
_{A_h\text{ term}}
\underalign{\text{(1)}}{\le}\delta(n+D_0).
\end{aligned}
\]
This is \eqref{eq:pred-13} in Proposition~\ref{prop:ricci-reward-master-inequality},
which is the starting point for (4) and (5) below.

\item Apply (3) to maximizing discs for $h=h_0$,
compare them with arbitrary finite-area holomorphic discs, and let $\delta\downarrow0$. This gives
\[
P_{\operatorname{Ric}(h_0)}(f)\le C
\]
with an upper bound independent of the disc $f$
(see Corollary~\ref{cor:uniform-ricci-poisson-bound}).
We then use discs that are holomorphic on a neighborhood of the closed disc and set
\[
u(x):=\inf_{g(0)=x}\bigl\{-P_{\operatorname{Ric}(h_0)}(g)\bigr\}
\]
Magn\'usson's disc envelope theorem \cite[Theorem~1.1]{Ch7Magnusson} gives
\(
-\operatorname{Ric}(h_0)
+\sqrt{-1}\,\partial\bar\partial u\ge0
\)
and $-C\le u\le0$. So $K_X$ admits a singular metric whose Lelong numbers vanish everywhere.
Hence $K_X$ is nef
(see Theorem~\ref{thm:canonical-nef-envelope}).

\item Fix $B>0$ such that $R_h\ge-B$.
For the functional $J_\delta(f):=
d(f)
+P_q(f)
-\delta A_h(f)$,
let $f_\delta $ be a maximizing disc. Combining the lower bounds in (1), we obtain
\[
n\exp\!\left(
\frac{-D_0-\log\frac{\omega_h^n}{\omega_{h_0}^n}(x)}n
\right)
\underalign{\text{(1)}}{\le}
n\exp\!\left(
\frac{P_{\operatorname{Ric}(h)}(f_\delta)
-\log\frac{\omega_h^n}{\omega_{h_0}^n}(x)}n
\right)
\underalign{\text{(3) \& }R_h\ge-B}{\le}
B+\delta(n+D_0).
\]
Since the left-hand side is independent of $\delta$, letting $\delta\downarrow0$ gives
\[
\frac{\omega_h^n}{\omega_{h_0}^n}(x)\ge
e^{-D_0}\left(\frac nB\right)^n
\]
This is Theorem~\ref{thm:scalar-volume-by-finite-prediction}.

\item Using the nefness in (4) and the Monge--Amp\`ere equation,
we construct, for each $\varepsilon>0$, a K\"ahler metric $h_\varepsilon$ satisfying
\[
[\omega_{h_\varepsilon}]
=2\pi c_1(K_X)+\varepsilon[\omega_{h_0}],
\qquad
\operatorname{Ric}(h_\varepsilon)
=-\omega_{h_\varepsilon}+\varepsilon\omega_{h_0}
\]
as in
\cite[Proposition~8(i)]{Ch8WuYau}.
Then
\[
R_{h_\varepsilon}
=-n+\underbrace{\varepsilon
\operatorname{tr}_{h_\varepsilon}h_0}_{\ge0}
\ge-n
\quad\Longrightarrow\quad
\omega_{h_\varepsilon}^n
\underalign{\text{apply (5) with }B=n}{\ge}
e^{-D_0}\omega_{h_0}^n.
\]
Integrating and letting $\varepsilon\downarrow0$, we obtain
\[
(2\pi)^nc_1(K_X)^n
=\lim_{\varepsilon\downarrow0}\int_X\omega_{h_\varepsilon}^n
\ge e^{-D_0}\int_X\omega_{h_0}^n>0.
\]
The nefness in (4) and this positivity imply that $K_X$ is big
\cite[Corollary~6.19]{DemaillyAnalyticMethods}
(Theorem~\ref{thm:canonical-big-volume}).
\end{enumerate}

Finally, the bigness of $K_X$ and the K\"ahler assumption imply that $X$ is projective,
and Brody hyperbolicity excludes rational curves.
Thus the standard argument for projective manifolds
with nef and big canonical bundle and no rational curves
shows that $K_X$ is ample
(cf.~ \cite[p.~219]{DebarreHigherDimensional}).

For (4) and (6), there is an alternative proof of nefness and bigness using the K\"ahler--Ricci flow, as in Nomura~\cite[Theorems~1.1--1.2]{Nomura}. We present this in Appendix~\ref{sec:krf-appendix}.
As a remark, surprisingly, this argument (2) is also similar to the analytic construction of rational curves
on Fano manifolds in \cite{DuGuoXie2026}.

\subsection{Comparison with the OpenAI manuscript}
Some parts of the original argument were unclear to the author. We therefore modified these arguments with the assistance of ChatGPT 6 Astra.
The following table summarizes the main differences between the OpenAI manuscript~\cite{OpenAI2026Canonical} and this paper.

\begin{longtable}{p{0.14\linewidth}p{0.38\linewidth}p{0.38\linewidth}}
\hline
Topic & OpenAI manuscript & This note\\
\hline
\endfirsthead
\hline
Topic & OpenAI manuscript & This note\\
\hline
\endhead

Holo\-morphic frames
&
Uses local $\bar\partial$ equations and the Oka--Grauert principle
\cite[Lemma~4.1]{OpenAI2026Canonical}.
&
Uses a right inverse fixing the value at one point and induction on the rank
(Theorem~\ref{thm:disc-frame-variation-hessian}(1)).
\\[1ex]

Boundary orthonormalization
&
Uses matrix factorization via a weighted Hardy space
to obtain $W^{1,p_0}$ regularity for some $2<p_0\le p$
\cite[Proposition~4.2]{OpenAI2026Canonical}.
&
Uses Lempert's factorization theorem and reflection across the circle
to obtain $W^{1,p}$ regularity
(Theorem~\ref{thm:disc-frame-variation-hessian}(2)).
\\[1ex]

Defor\-mations fixing the center
&
Applies the complex Banach implicit function theorem to a nonlinear $\bar\partial$ equation
\cite[Proposition~4.3]{OpenAI2026Canonical}.
&
Uses the complex Banach manifold structure on the space of holomorphic maps and evaluation at the center
(Theorem~\ref{thm:disc-center-fixed-deformation}).
\\[1ex]

Second variation
&
States trace identities obtained by summing over the deformation directions
\cite[Proposition~5.2]{OpenAI2026Canonical}.
&
Proves the second variation formula in each direction and then takes the sum
(Theorem~\ref{thm:disc-second-variation},
Theorem~\ref{thm:signed-disc-determinant-trace}).
\\[1ex]

Boundary integrals
&
Along a suitable sequence of radii,
compares the boundary derivatives of $\operatorname{tr}G$ and $\log\det G$
\cite[Proof of Proposition~5.2]{OpenAI2026Canonical}.
&
Shows that the circular mean of $\operatorname{tr}G-\log\det G-n$ is $o(1-r)$,
and hence that the integral of its Laplacian vanishes
(Claim~\ref{claim:ch5-curvature-and-center}).
\\[1ex]

Inequality at a maximum
&
Obtains a trace inequality for maximizing discs
\cite[Corollary~5.3]{OpenAI2026Canonical}.
&
Uses the arithmetic--geometric mean inequality for locally maximizing discs
to obtain an inequality involving a determinant
(Theorem~\ref{thm:signed-disc-determinant-trace}).
\\[1ex]

Volume lower bound
&
Estimates the volume form from below under the assumption $\operatorname{Ric}(h)\ge-\omega_h$
\cite[Proposition~7.1]{OpenAI2026Canonical}.
&
Estimates the volume form from below using a scalar curvature lower bound $R_h\ge-B$
(Theorem~\ref{thm:scalar-volume-by-finite-prediction}).
\\[1ex]

Nefness
&
On the total space of $\det TX$ with the zero section removed,
uses the Poletsky--Rosay theorem to construct a bounded potential
\cite[Propositions~6.2, 6.4]{OpenAI2026Canonical}.
&
Applies Magn\'usson's disc envelope theorem on $X$
to construct a bounded potential
(Theorem~\ref{thm:canonical-nef-envelope}).
\\[1ex]

Top self-intersection and bigness
&
Uses the Monge--Amp\`ere equation and a volume lower bound to prove positivity of the top self-intersection,
and combines this with nefness to obtain bigness
\cite[Proposition~8.2, Lemma~9.1]{OpenAI2026Canonical}.
&
Uses Wu--Yau's construction of metrics and the volume lower bound in this paper
to obtain bigness by the same strategy
(Theorem~\ref{thm:canonical-big-volume}).
\\[1ex]

K\"ahler--Ricci flow
&
Not used in this manuscript
\cite{OpenAI2026Canonical}.
&
Proves long-time existence and obtains bigness from nefness and positive top self-intersection
(Theorem~\ref{thm:krf-nef-big}).
\\
\hline
\end{longtable}

\medskip
\noindent{\bfseries Acknowledgments.}
M.\,I.\ would like to thank Professors Yusaku Tiba, Katsutoshi Yamanoi, and Ryoichi Kobayashi for pointing out the connections between the functionals \(P_q\) and \(A_h\) and Nevanlinna theory.
He  also would like to thank Shiyu Zhang for sharing his recent preprint \cite{Zhang2026VolumeSchwarz}.

M.\,I.\ was supported by the Grant-in-Aid for Early Career Scientists, No.\ 22K13907. 
He also thanks OpenAI for providing complimentary access to ChatGPT
through the ChatGPT for Academic Researchers program.

\medskip
\noindent{\bfseries Use of AI.}
The author first asked ChatGPT 6 Astra to explain the paper \cite{OpenAI2026Canonical}.
The explanation made by ChatGPT 6 Astra was very difficult to read, and the author could not understand it at all.
In the end, the author had to work through each theorem and proposition one by one.
% and he wasted a whole week. 

Nevertheless, ChatGPT 6 Astra helped the author understand the proofs
and write out each step carefully and in detail.
It also helped the author identify the connection between \cite{OpenAI2026Canonical} and \cite{DuGuoXie2026}.
%At present, AI seems to be good at proving results but poor at writing papers.Proofs written by AI are very difficult to read,so it seems that humans still need to write explanations.

The author used ChatGPT 6 Astra 
to improve the writing, search the literature, and discuss
and check mathematical arguments and calculations.
The author reviewed all AI-assisted material included in
the paper and checked the arguments, calculations, and
references himself.
The author takes full responsibility for the results and
proofs in this paper and for citing previous work accurately
and fully.

\section{Preliminaries}
We denote  by $\D$ the unit disc.
Unless otherwise stated, $X$ is an $n$-dimensional compact K\"ahler manifold and $h_0$ is a fixed K\"ahler metric.
We write $\omega_h$ for the K\"ahler form associated with a K\"ahler metric $h$. Throughout this paper, $TX$ denotes the holomorphic tangent bundle.
We use $\Ric(h_0) \in 2\pi c_1(X)$.
\subsection{Sobolev spaces}
\begin{defn}

Let
\(
\mathbb D=\{z\in\mathbb C:|z|<1\}
\)
be the unit disc, and write $z=x+\sqrt{-1}y$ and $dS=dx\,dy$.
We consider complex-valued functions and assume $2<p<\infty$.

The $L^p$-space $L^p(\mathbb D)$ consists of measurable functions $b : \D \to \C$ satisfying
\[
\|b\|_{L^p(\mathbb D)}
:=
\left(\int_{\mathbb D}|b|^p\,dS\right)^{1/p}
<\infty
\]
The Sobolev space $W^{1,p}(\mathbb D)$ consists of functions whose first weak derivatives, as well as the functions themselves,
belong to $L^p$; that is,
\[
W^{1,p}(\mathbb D)
=
\left\{
u\in L^p(\mathbb D):
\partial_xu,\partial_yu\in L^p(\mathbb D)
\right\}
\]
We equip it with the norm
\[
\|u\|_{W^{1,p}}
=
\|u\|_{L^p}
+\|\partial_xu\|_{L^p}
+\|\partial_yu\|_{L^p}
\]

\end{defn}

\begin{defn}
\label{def:sobolev-maps-into-X}
Let $X$ be a smooth compact complex manifold of complex dimension $n$,
and let $2<p<\infty$.
We say that a measurable map $f:\mathbb D\to X$ belongs to
$W^{1,p}(\mathbb D,X)$ if it has a continuous representative
extending to $\overline{\mathbb D}$,
and there exist finitely many open discs $B_i\subset\mathbb C$
and holomorphic coordinate charts
$\varphi_i:U_i\to\varphi_i(U_i)\subset\mathbb C^n$
such that:
\begin{itemize}
\item
The discs cover $\overline{\mathbb D}$, 
\(
\overline{\mathbb D}\subset\bigcup_{i=1}^N B_i\), and
\(f\bigl(\overline{\mathbb D}\cap\overline{B_i}\bigr)
\subset U_i.
\)

\item
Each coordinate representation satisfies
\[
\varphi_i\circ f
\in W^{1,p}(\mathbb D\cap B_i,\mathbb C^n).
\]
Equivalently, writing $z=x+\sqrt{-1}y$, for each
$\alpha=1,\ldots,n$,
\[
\varphi_i^\alpha\circ f,\qquad
\partial_x(\varphi_i^\alpha\circ f),\qquad
\partial_y(\varphi_i^\alpha\circ f)
\in L^p(\mathbb D\cap B_i),
\]
where the derivatives are weak derivatives.
\end{itemize}
\end{defn}
Since $p>2$, every map in $W^{1,p}(\mathbb D,X)$ has a unique continuous representative on $\overline{\mathbb D}$ by the Sobolev embedding theorem \cite[Theorem~3.49(3)]{Ch4Hunter}, so continuity need not be imposed as a separate assumption.

\subsection{Green's Formula}

\begin{lem}
\label{lem:green-normalization}
For a smooth function $u : \D \to \R$ and $0<r_0<1$, we have
\begin{equation}
\label{eq-defn-Mv}
\begin{aligned}
\frac1{2\pi}\int_0^{2\pi}u(r_0e^{\sqrt{-1}\theta})\,d\theta-u(0)
&=\frac2\pi\int_{|z|<r_0}\log\frac{r_0}{|z|}\,u_{z\bar z}\,dS,\\
\frac1\pi\int_{|z|<r_0}u_{z\bar z}\,dS
&=\frac{r_0}{4\pi}\int_0^{2\pi}u_r(r_0e^{\sqrt{-1}\theta})\,d\theta.
\end{aligned}
\end{equation}

Moreover, the following holds.
\begin{itemize}
\item if $u\in C^0(\overline{\mathbb D})$ and
$u_{z\bar z}\in L^a(\mathbb D)$ for $a>1$, then
\begin{equation}
\label{eq-lem:green-normalization-2}
\frac1{2\pi}\int_0^{2\pi}u(e^{\sqrt{-1}\theta})\,d\theta-u(0)
=\frac2\pi\int_{\mathbb D}\log\frac1{|z|}\,u_{z\bar z}\,dS.
\end{equation}
\item 
if $\Delta u\in L^1(\mathbb D)$, we have
\begin{equation}
\label{eq:pred-3}
\frac1\pi\int_{\mathbb D}u_{z\bar z}\,dS
%=\frac1{4\pi}\int_{\mathbb D}\Delta u\,dS
=\lim_{r_0\uparrow1}\frac{r_0}{4\pi}
\left.\frac d{dr}\int_0^{2\pi}u(re^{\sqrt{-1}\theta})\,d\theta\right|_{r=r_0}.
\end{equation}
\end{itemize}
\end{lem}

\begin{proof}
We use the Laplacian symbol  \(\Delta u :=4u_{z\bar z}\).
In \cite[Ch.I, 4.B]{Ch4Demailly},
\(d\lambda=dS\) is the volume element and \(d\sigma\) is the boundary area element. In particular, on $|z|=r_0$,
we have \(d\sigma=r_0 d \theta\).
The Green and Poisson kernels in \cite[Ch.I, (4.2), (4.6)]{Ch4Demailly} are
\begin{equation}
\label{eq-defn-gree-poison}
G_{r_0}(0,z)=\frac1{4\pi}\log\frac{|z|^2}{r_0^2}
\underalign{\cite[(4.2)]{Ch4Demailly}}{=}\frac1{2\pi}\log\frac{|z|}{r_0},
\quad
P_{r_0}(0,z)
\underalign{\cite[(4.6)]{Ch4Demailly}}{=}
\frac1{2\pi |z|}.
\end{equation}
Hence
Demailly's Green--Riesz formula in \cite[Ch.I. 4.B]{Ch4Demailly} gives
\[
\begin{aligned}
u(0)
&\underalign{\cite[(4.5)]{Ch4Demailly}}{=}\frac1{2\pi}\int_{|z|<r_0}
    \underbrace{\log\frac{|z|}{r_0}}_{=2\pi G_{r_0}(0,z) \text{ by \eqref{eq-defn-gree-poison}}}\,\Delta u\,dS
+\int_0^{2\pi}u(r_0e^{\sqrt{-1}\theta})
\underbrace{\frac1{2\pi r_0}}_{=P_{r_0}(0,z) \text{ by \eqref{eq-defn-gree-poison}}}
r_0\,d\theta
\end{aligned}
\]
Also, Green's formula~\cite[(4.4)]{Ch4Demailly} gives
\[
4\int_{|z|<r_0}u_{z\bar z}\,dS
=
\int_{|z|<r_0} \Delta u\,dS
\underalign{\text{\cite[(4.4)]{Ch4Demailly}}}{=}
r_0\int_0^{2\pi}
u_r(r_0e^{\sqrt{-1}\theta})
\,d\theta.
\]
This proves both identities in \eqref{eq-defn-Mv}.

If $u_{z\bar z}\in L^a$ with $a>1$,
then $\log(1/|z|)\in L^{a/(a-1)}$. Thus H\"older's inequality shows that
$\log\frac1{|z|}\,|u_{z\bar z}|$ is
integrable.
By the dominated convergence theorem,
letting $r_0\uparrow1$ gives \eqref{eq-lem:green-normalization-2}.
Similarly, if $\Delta u\in L^1$,
letting $r_0\uparrow1$ in the second identity of \eqref{eq-defn-Mv}
gives \eqref{eq:pred-3}.
\end{proof}

\begin{cor}
For any $u\in W^{1,p}(\mathbb D)$ satisfying $u_{z\bar z}\in L^{p/2}(\mathbb D)$
with $p>2$, we obtain 
\begin{equation}
\label{eq:green-radial-product}
\frac1\pi\int_{\mathbb D}\partial_z\partial_{\bar z}(|z|^2u)\,dS
=\frac1{2\pi}\int_0^{2\pi}u(e^{\sqrt{-1}\theta})\,d\theta
+\frac1\pi\int_{\mathbb D}u_{z\bar z}\,dS.
\end{equation}
\end{cor}
\begin{proof}
First, if $u$ is smooth, then
\begin{align*}
\frac1\pi\int_{|z|<r_0}
\partial_z\partial_{\bar z}(|z|^2u)\,dS
&\underalign{\eqref{eq-defn-Mv}}{=}
\frac{r_0}{4\pi}\left.\frac d{dr}
\left(r^2\int_0^{2\pi}u(re^{\sqrt{-1}\theta})\,d\theta\right)\right|_{r=r_0}\\
&\underalign{\text{}}{=}
\frac{r_0^2}{2\pi}\int_0^{2\pi}u(r_0e^{\sqrt{-1}\theta})\,d\theta
+r_0^2\underbrace{
\frac{r_0}{4\pi}\left.\frac d{dr}\int_0^{2\pi}u(re^{\sqrt{-1}\theta})\,d\theta\right|_{r=r_0}
}_{
=\frac1\pi\int_{|z|<r_0}u_{z\bar z}\,dS}.
\end{align*}
By the dominated convergence theorem,
letting $r_0\uparrow1$ gives \eqref{eq:green-radial-product}.

Fix $r_0<1$ and let $\rho_\varepsilon$ be a mollifier.
Define $u_\varepsilon:=\rho_\varepsilon*u$ on a neighborhood of $\{|z|\le r_0\}$.
As $\varepsilon\downarrow0$, we have
\[
u_\varepsilon\to u\ \text{uniformly},\qquad
\nabla u_\varepsilon\to\nabla u\ \text{in }L^p,\qquad
\Delta u_\varepsilon\to\Delta u\ \text{in }L^{p/2}.
\]
Thus $\Delta(|z|^2u_\varepsilon)\to\Delta(|z|^2u) \in L^1$.
Thus, letting $\varepsilon\downarrow0$ in the identity for smooth functions gives
\[
\frac1\pi\int_{|z|<r_0}(|z|^2u)_{z\bar z}\,dS
=
\frac{r_0^2}{2\pi}\int_0^{2\pi}u(r_0e^{\sqrt{-1}\theta})\,d\theta
+\frac{r_0^2}{\pi}\int_{|z|<r_0}u_{z\bar z}\,dS.
\]
Finally, using the continuity of $u$ at the boundary and $(|z|^2u)_{z\bar z} \in L^1$
%\[
%\Delta(|z|^2u)
%=\underbrace{|z|^2\Delta u+4(xu_x+yu_y)+4u}_{\in L^1(\mathbb D)}
%\]
we let $r_0\uparrow1$ to obtain \eqref{eq:green-radial-product}.
\end{proof}

\section{Existence of maximizing discs}

\begin{defn}
\label{defn-APJ-funcional}
Let \(X\) be a compact K\"ahler manifold and $h_0$ a fixed metric.
For a holomorphic disc \(f:\mathbb D\to X\) and a real \((1,1)\)-form \(q\), we define the following functionals:
\begin{itemize}
\item Set $d(f):=|f'(0)|_{h_0}^2$.
\item For $q=\sqrt{-1}\sum_{a,b=1}^n q_{a\bar b}\,dw^a\wedge d\bar w^b$,
we write $q(v,v):=\sum_{a,b=1}^n q_{a\bar b}v^a\overline{v^b}$. We set
$$
A_q(f):=\frac1\pi\int_{\mathbb D}q(f',f')\,dS
%=\frac1\pi\int_{\mathbb D}\sum_{a,b=1}^nq_{a\bar b}(f(z))\,(f^a)'(z)\,\overline{(f^b)'(z)}
$$
In particular, if $q=\omega_{h_0}$ (that is, if $q$ is the $(1,1)$-form induced by $h_0$), then $q(f',f')=h_0(f',f')=|f'(z)|_{h_0}^2$, so
$$
A_{h_0}(f)=\frac1\pi\int_{\mathbb D}|f'(z)|_{h_0}^2\,dx\,dy
=
\frac1{2\pi}\int_{\mathbb D}f^*\omega_{h_0}.
$$
We say that $f$ has finite area if $A_{h_0}(f)<+\infty$.
\item We define $P_q(f)$ by 
$$
P_q(f):=\frac2\pi\int_{\mathbb D}\log\frac1{|z|}q(f',f')\,dS
$$
In particular, if $q=\omega_{h_0}$, then
$$
P_{h_0}(f)=\frac2\pi\int_{\mathbb D}
\log\frac1{|z|}
|f'(z)|_{h_0}^2\,dS
$$
%\item \(A_{q,r}(f)=\pi^{-1}\int_{|z|<r}q(f',f')\,dS\)とおく
%\item 
%\begin{equation}
%\label{eq:pred-1}
%J(f)=d(f)+P_q(f)-\delta A_h(f),
%\qquad \delta>0.
%\tag{2.1}
%\end{equation}

\end{itemize}

\end{defn}
In this notation, $\omega_h(v,v)=h(v,v)$ and
we write $A_h:=A_{\omega_h}$ and $P_h:=P_{\omega_h}$.

\begin{lem}[{\cite[Lemma~3.1]{OpenAI2026Canonical}, cf. \cite[Corollary~1.8]{Demailly2012Hyperbolic}}]
\label{lem:brody-bound}
Assume that \(X\) is Brody hyperbolic, and let $h_0$ be a fixed K\"ahler metric.
Then there exists $D_0 >0$ such that every holomorphic disc satisfies
\[d(f)\le D_0
\quad
|f'(z)|_{h_0}\le\sqrt {D_0}/(1-|z|)
\]

Moreover, every sequence of holomorphic discs $f_j:\mathbb D\to X$ has a
subsequence converging locally uniformly to a holomorphic map $f$,
with derivatives converging locally uniformly to $f'$
\end{lem}

\begin{proof}
If $d(f)$ were unbounded,
there would exist a sequence of holomorphic discs with $|f_j'(0)|_{h_0}\to\infty$.
The consequence of Brody's reparametrization lemma stated in
\cite[Corollary~1.8]{Demailly2012Hyperbolic} would then yield
a nonconstant holomorphic map $g:\mathbb C\to X$, a contradiction.
Therefore, there exists $D_0>0$ such that
\begin{equation}
\label{eq:brody-proof-center-bound}
d(f)=|f'(0)|_{h_0}^2\le D_0
\end{equation}
for every holomorphic disc.
For any $z\in\mathbb D$,
the map $F_z(w):=f(z+(1-|z|)w)$ is also a holomorphic disc on $\mathbb D$. Hence
\[
(1-|z|)^2|f'(z)|_{h_0}^2
\underalign{\text{chain rule}}{=}
\underbrace{|F_z'(0)|_{h_0}^2}_{=d(F_z)}
\underalign{\eqref{eq:brody-proof-center-bound}}{\le}
D_0.
\]
Thus $|f'(z)|_{h_0}\le\sqrt{D_0}/(1-|z|)$.
For the last statement, the preceding estimate and compactness of $X$ give
$f_{j_k}\to f$ locally uniformly by Arzel\`a--Ascoli and a diagonal argument.
Cauchy's integral formula in local holomorphic coordinates shows that
$f$ is holomorphic and $f'_{j_k}\to f'$ locally uniformly.
\end{proof}

\begin{prop}[{\cite[Propositions~3.2, 3.4]{OpenAI2026Canonical}}]
\label{prop:disc-maximizer}
Assume that \(X\) is Brody hyperbolic, and let $h_0$ be a fixed K\"ahler metric.
Fix a smooth K\"ahler metric $h$, a smooth closed real $(1,1)$-form $q$,
$\delta>0$, and $x\in X$, and define
\[
\mathcal A_x
:=
\left\{
f:\mathbb D\to X\ \text{holomorphic}:
f(0)=x,\,A_{h_0}(f)<\infty
\right\},
\quad
J_\delta(f):=d(f)+P_q(f)-\delta A_h(f)
\]
Then the following hold.
\begin{enumerate}
\item
There exists $f_*\in\mathcal A_x$ such that
\[
0\le J_\delta(f_*)=\sup_{f\in\mathcal A_x}J_\delta(f)<\infty
\]
A disc $f_{*} :\D \to X$ attaining this maximum is called
a maximizing disc for $J_\delta$.

\item
Every maximizing disc $f_*$ extends continuously to the closed disc $\overline{\D}$,
and for every $2<p<4$ it satisfies
\[
f_*\in W^{1,p}(\mathbb D,X),
\quad \text{and} \quad 
\int_{\mathbb D}|f_*'(z)|_{h_0}^p\,dS(z)<\infty
\]
\end{enumerate}
\end{prop}

\begin{proof}
(1). Notice that $C$ denotes a constant independent of $f$ that may change from line to line. (Thus the constants $C$ in the inequalities need not be the same.)

Since $X$ is compact, the metrics $h$ and $h_0$ are uniformly comparable,
and there exists $C_q>0$ such that $|q(v,v)|\le C_qh(v,v)$.
Choose $r_0\in[1/2,1)$ such that $2C_q\log(1/r_0)\le\delta/2$.
Then we can take a constant $C>0$ such that
\begin{equation}
\label{eq:max-proof-coercivity}
\begin{aligned}
J_\delta(f)
&\underalign{}{=}
\underbrace{d(f)}_{\le C \text{ by Lemma~\ref{lem:brody-bound}}}+\frac1\pi
\underbrace{\int_{|z|<r_0}\left(2\log\frac1{|z|}\,q(f',f')-\delta h(f',f')\right)dS}_{ \le C \text{ by Lemma~\ref{lem:brody-bound} and $\log\frac1{|z|} \in L^{1}(\{ |z|<r_0\})$}}
 \\
&
\underalign{}{+}
\frac1\pi\int_{r_0<|z|<1}
\underbrace{\left(2\log\frac1{|z|}\,q(f',f') - \delta h(f',f')\right)}
_{\le-\frac{\delta}{2}h(f',f') \le 0}
dS
\end{aligned}
\end{equation}
Since $J_\delta=0$ for constant discs, we have $0\le\sup_{\mathcal A_x}J_\delta<\infty$.

Choose a maximizing sequence $f_j$ with $J_\delta(f_j)\ge-1$.
By Lemma~\ref{lem:brody-bound}, after passing to a subsequence,
$f_j$ converges locally uniformly to a holomorphic map $f_{*}:\mathbb D\to X$,
and the derivatives converge locally uniformly.

We show that $f_{*} \in \mathcal{A}_x$. The convergence gives $f_{*}(0)=x$. To prove finite area, we use
\eqref{eq:max-proof-coercivity} to obtain
\[
-1\le J_\delta(f_j)
\underalign{\eqref{eq:max-proof-coercivity} }{\le} C-\frac{\delta}{2\pi}
\int_{r_0<|z|<1}h(f_j',f_j')\,dS.
\]
Therefore,
\[
A_h(f_j)
=
\underbrace{\frac1\pi\int_{|z|<r_0}h(f_j',f_j')\,dS}
_{\le C\text{ by }h\le C_hh_0\text{ and Lemma~\ref{lem:brody-bound}}}
+\frac1\pi\int_{r_0<|z|<1}h(f_j',f_j')\,dS
\le C+\frac{2(C+1)}{\delta}.
\]
The right-hand side is independent of $j$, so $\sup_jA_h(f_j)<\infty$.
Fatou's lemma now gives $A_h(f_{*}) \le \liminf_jA_h(f_j)<\infty$. 
This shows $f_{*}\in\mathcal A_x$.

Finally, we show that $f_{*}$ attains the maximum of $J_\delta$.
The first two terms in \eqref{eq:max-proof-coercivity}
converge to their values at $f_{*}$ as $j \to \infty$.
For the third term in \eqref{eq:max-proof-coercivity}, we have
\begin{equation}
\label{eq:max-proof-coercivity-2} 
\begin{aligned}
&\frac1\pi\int_{r_0<|z|<1}
\left(2\log\frac1{|z|}\,q(f'_{*},f'_{*})-\delta h(f'_{*},f'_{*})\right)dS
\\
&\underalign{\text{Fatou's lemma}}{\ge}
\limsup_{j\to\infty}\frac1\pi\int_{r_0<|z|<1}
\left(2\log\frac1{|z|}\,q(f_j',f_j')-\delta h(f_j',f_j')\right)
dS.
\end{aligned}
\end{equation}
Hence
\[
\sup_{\mathcal A_x}J_\delta
\underalign{\text{ by $f_{*}\in\mathcal A_x$ }}{\ge}
J_\delta(f_{*})
\underalign{ \eqref{eq:max-proof-coercivity} \& \eqref{eq:max-proof-coercivity-2} }
{\ge}
\lim_{j\to\infty}J_\delta(f_j)
\underalign{\text{definition}}{=}
\sup_{\mathcal A_x}J_\delta.
\]
Thus $J_\delta(f_{*})=\sup_{\mathcal A_x}J_\delta$, proving (1).

(2). Let $f$ be a maximizing disc. Notice that $C$ may depend on $f$ and may change from line to line. 
Set $f_r(z):=f(rz)$. A change of variables and Fubini's theorem give
$d(f_r)=r^2d(f)$ and
\begin{equation}
\label{eq:max-proof-scaling}
\begin{aligned}
A_h(f_r)=\frac1\pi\int_{|z|<r}h(f',f')\,dS,
\quad
P_q(f_r)
=\frac2\pi\int_{|z|<r}
\underbrace{\log\frac r{|z|}}_{=\int_{|z|}^r ds/s}
q(f',f')\,dS
\underalign{\text{Fubini's theorem}}{=}
2\int_0^r A_q(f_s)\,\frac{ds}{s}.
\end{aligned}
\end{equation}
Therefore, for $1/2\le r<1$,
\begin{equation}
\label{eq:disc-maximizer-1}
\begin{aligned}
0\le\delta\bigl(A_h(f)-A_h(f_r)\bigr)
&\underalign{J_\delta(f)\ge J_\delta(f_r)\text{ \& }\eqref{eq:max-proof-scaling}}{\le}
(1-r^2)d(f)+2\int_r^1A_q(f_s)\,\frac{ds}{s}\\
&\underalign{|q(v,v)|\le C_qh(v,v) \&  A_h(f_s)\le A_h(f)}{\le}
(1-r^2)d(f)
+2C_qA_h(f) \log(1/r)
\le C(1-r).
\end{aligned}
\end{equation}
Comparing $h$ and $h_0$, we obtain
\begin{equation}
\label{eq:max-proof-area-tail}
\int_{r<|z|<1}|f'(z)|_{h_0}^2\,dS
\underalign{}{\le}
C\int_{r<|z|<1}h(f',f')\,dS
\underalign{\eqref{eq:max-proof-scaling}}{=}
C\left(A_h(f)-A_h(f_r)\right)
\underalign{\eqref{eq:disc-maximizer-1}}{\le}
C(1-r)
\end{equation}
To define the boundary values by $\overline f(\zeta):=\lim_{r\uparrow1}f(r\zeta)$,
we show that this limit exists uniformly
for $\zeta\in\partial\mathbb D$.

First, we estimate the derivative at $z\in\mathbb D$ with $|z|\ge3/4$.
By Lemma~\ref{lem:brody-bound},
we may fix $0<c<1/2$ sufficiently small so that, for every $z\in\mathbb D$,
the image $f(B(z,c(1-|z|)))$ lies in a single coordinate neighborhood.
In particular, $h_0$ and the Euclidean metric can be compared uniformly.
Using the mean value inequality~\cite[Chapter~I, (4.12)(c)]{Ch4Demailly} and $B(z,c(1-|z|))\subset
\{|w|>1-(1+c)(1-|z|)\}$, we obtain
\begin{equation}
\label{eq:max-proof-derivative}
\begin{aligned}
|f'(z)|_{h_0}^2
&\underalign{\text{ \cite[Chapter~I, (4.12)(c)]{Ch4Demailly}}}{\le}
\frac{C}{(1-|z|)^2}
\int_{B(z,c(1-|z|))}|f'|_{h_0}^2\,dS\\
&\underalign{
}{\le}
\frac{C}{(1-|z|)^2}
\int_{1-(1+c)(1-|z|)<|w|<1}|f'(w)|_{h_0}^2\,dS
\underalign{\eqref{eq:max-proof-area-tail}}{\le}
\frac{C}{1-|z|}.
\end{aligned}
\end{equation}
Therefore, for $|\zeta|=1$ and $3/4\le r<s<1$,
\[
\operatorname{dist}_{h_0}(f(r\zeta),f(s\zeta))
\underalign{\eqref{eq:max-proof-derivative}}{\le}
C\int_r^s\frac{dt}{\sqrt{1-t}}
\le C\sqrt{1-r}.
\]
The completeness of $X$ implies uniform existence of the limit, so the boundary values are defined by $\overline f(\zeta):=\lim_{r\uparrow1}f(r\zeta)$.

We show that $\overline f$ is continuous on the closed disc.
If $r\zeta\to\zeta_0\in\partial\mathbb D$ with $|\zeta|=1$,
then $r\to1$ and
\[
|\zeta-\zeta_0|
\le|\zeta-r\zeta|+|r\zeta-\zeta_0|
=(1-r)+|r\zeta-\zeta_0|\longrightarrow0.
\]
Therefore,
\begin{align*}
\operatorname{dist}_{h_0}(f(r\zeta),\overline f(\zeta_0))
&\underalign{\text{}}{\le}
\operatorname{dist}_{h_0}(f(r\zeta),\overline f(\zeta))
+\operatorname{dist}_{h_0}(\overline f(\zeta),\overline f(\zeta_0))\\
&\underalign{\text{}}{\le}
\underbrace{C\sqrt{1-r}}_{\to0}
+\underbrace{
\operatorname{dist}_{h_0}(\overline f(\zeta),\overline f(\zeta_0))
}_{\to0\text{ by continuity of the boundary values}}
\longrightarrow0.
\end{align*}
Thus $\overline f$ is continuous on the closed disc.

Finally, set $E_j:=\{2^{-j-1}<1-|z|<2^{-j}\}$ for $j\ge2$.
For $2<p<4$, we have
\begin{align*}
\int_{3/4<|z|<1}|f'|_{h_0}^p\,dS
&\underalign{}{\le}
\sum_{j\ge2}
\underbrace{\left(\sup_{E_j}|f'|_{h_0}\right)^{p-2}}
_{ \le (C 2^{j+1})^{\frac{p-2}{2}} \text{ by \eqref{eq:max-proof-derivative}}}
\cdot
\underbrace{\int_{E_j}|f'|_{h_0}^2\,dS}_{\le C2^{-j}\text{ by \eqref{eq:max-proof-area-tail}}}
%\\
%&\underalign{\eqref{eq:max-proof-derivative}}{\le}
%C_{f,p}\sum_{j\ge2}
%\underbrace{2^{j(p-2)/2}}_{\text{derivative bound}}
%\underbrace{2^{-j}}_{\text{area bound}}
%=
\le 
C
\underbrace{\sum_{j\ge2}2^{-\frac{j(4-p)}{2}}}_{\text{converges since }p<4}
<\infty.
\end{align*}

Since the derivative is bounded away from the boundary,
we have $\int_{\mathbb D}|f'|_{h_0}^p\,dS<\infty$.
By continuity on the closed disc, we may use finitely many coordinate charts on the target.
Each coordinate representation is bounded, and its first weak derivatives belong to $L^p$.
Hence $f\in W^{1,p}(\mathbb D,X)$, proving (2).
\end{proof}

\begin{cor}[{\cite[Lemma~3.3, Proposition~3.2]{OpenAI2026Canonical}}]
\label{lem:signed-reward-source-bounds}
Assume that \(X\) is Brody hyperbolic, and let $h_0$ be a fixed K\"ahler metric.
Let $D_0$ be the constant in Lemma~\ref{lem:brody-bound}.

With the notation of Proposition~\ref{prop:disc-maximizer},
every maximizing disc $f\in\mathcal A_x$ satisfies
\begin{equation}
\label{eq:pred-4}
\begin{aligned}
A_q(f) \ge-d(f)\ge-D_0,
\qquad
P_q(f)\ge-d(f)+\delta A_h(f)\ge-D_0
\end{aligned}
\end{equation}

\end{cor}

\begin{proof}
For $0<s<1$, define $f_s:\mathbb D\to X$ by $f_s(z):=f(sz)$.
Dividing \eqref{eq:disc-maximizer-1} by $1-r^2$, we obtain
\begin{equation}
\label{eq:source-bound-dilation}
0
\underalign{\eqref{eq:disc-maximizer-1}}{\le}
d(f)+\frac{2}{1-r^2}\int_r^1A_q(f_s)\,\frac{ds}{s}.
\end{equation}
Since $f$ has finite area and $|q(v,v)|\le Ch_0(v,v)$,
the dominated convergence theorem gives $A_q(f_s)\to A_q(f)$. Therefore, for any $r < s < 1$, 
\begin{equation}
\label{eq:source-bound-average-limit}
\begin{aligned}
\left|
\frac1{1-r}\int_r^1\frac{A_q(f_s)}s\,ds-A_q(f)
\right|
&\underalign{\text{by $\frac1{1-r}\int_r^1ds=1$}}{=}
\left|
\frac1{1-r}\int_r^1
\left(\frac{A_q(f_s)}s-A_q(f)\right)ds
\right|\\
&\underalign{}{\le}
\underbrace{\frac1{1-r}\int_r^1ds}_{=1}
\cdot 
\underbrace{
\sup_{r<s<1}\left|\frac{A_q(f_s)}s-A_q(f)\right|
}_{\to0\text{ since }A_q(f_s)/s\to A_q(f)}
\underset{r \to 1}{\longrightarrow}0.
\end{aligned}
\end{equation}
Letting $r\uparrow1$ in \eqref{eq:source-bound-dilation}
and using \eqref{eq:source-bound-average-limit} and $2/(1+r)\to1$, we obtain
\[
A_q(f)
\underalign{
\eqref{eq:source-bound-dilation} \& 
\eqref{eq:source-bound-average-limit}
}{\ge}
-d(f)
\underalign{\text{Lemma~\ref{lem:brody-bound}}}{\ge}
-D_0.
\]
\[
P_q(f)
\underalign{\text{}}{=}
\underbrace{J_\delta(f)}_{\ge0 \text{ by Proposition~\ref{prop:disc-maximizer}(1)}}-d(f)+\underbrace{\delta A_h(f)}_{\ge0}
\underalign{\text{}}{\ge}
-d(f)
\underalign{\text{Lemma~\ref{lem:brody-bound}}}{\ge}-D_0.
\]
\end{proof}

\section{Holomorphic frames and deformations fixing the center}
\label{sec:frames-unified}

The following two theorems form the most difficult and technical part of this note.

\begin{thm}[{cf.~\cite[Lemma~4.1, Proposition~4.2]{OpenAI2026Canonical}}]
\label{thm:disc-frame-variation-hessian}
Let $X$ be a smooth compact complex manifold of complex dimension $n\ge1$,
and let $2<p<\infty$.
Let $E\to\overline{\mathbb D}$ be a complex vector bundle of rank $r\ge1$. 

Assume that for any sufficiently small open disc $U\subset\mathbb C$, there exist local frames $e = (e_1, \ldots, e_r)$ defined near
$\overline{\mathbb D}\cap\overline U$ such that each $e_a$ is holomorphic on $\mathbb D\cap U$,
continuous on $\overline{\mathbb D}\cap U$, and belongs to
$W^{1,p}(\mathbb D\cap U,E)$. 
Then the following hold
\begin{enumerate}
\item $E$ admits a frame $H=(H_1,\ldots,H_r)$ such that 
each $H_i:\overline{\mathbb D}\longrightarrow E$ is continuous and holomorphic on $\D$
and, for every $z\in\overline{\mathbb D}$, the map
\[
H(z):\mathbb C^{r}\longrightarrow E_z,
\qquad
(c_1,\ldots,c_r)\longmapsto\sum_{i=1}^{r} c_iH_i(z)
\]
is a complex linear isomorphism.
Moreover, in local frames, the matrix representations satisfy
\[
H,H^{-1}\in W^{1,p}(\mathbb D)
\cap C^0(\overline{\mathbb D})
\] 
(See the remark below for the precise meaning.)

\item Let $h$ be a continuous Hermitian metric on $E$, positive definite on
$\overline{\mathbb D}$, such that, in every local frame
$e=(e_1,\ldots,e_r)$ as above,
\[
h(e_a,e_b)\in
W^{1,p}(\mathbb D\cap U,\mathbb C)
\cap C^0(\overline{\mathbb D}\cap U,\mathbb C)
\qquad(1\le a,b\le r).
\]
Then there exists a frame $B=(B_1,\ldots,B_r)$ with the same regularity and invertibility as in (1),
 satisfying
\[
h_z\bigl(B_i(z),B_j(z)\bigr)=\delta_{ij}
\qquad
(z\in\partial\mathbb D,\ 1\le i,j\le r)
\]
In particular, its local matrix representations satisfy
\[
B,B^{-1}\in W^{1,p}(\mathbb D)
\cap C^0(\overline{\mathbb D})
\]

\end{enumerate}

\end{thm}
\begin{rem}
For a frame $H=(H_1,\ldots,H_r)$, using a local frame $e = (e_1, \ldots, e_r)$ as in Theorem~\ref{thm:disc-frame-variation-hessian}, write
\[
H_j(z)=\sum_{a=1}^r \widetilde H_{aj}(z)e_a(z)
\qquad
(z\in\overline{\mathbb D}\cap U,\ 1\le j\le r).
\]
The notation
"$H\in W^{1,p}(\mathbb D)\cap C^0(\overline{\mathbb D})$"
means that, in every such local frame,
\[
\widetilde H_{aj}
\in W^{1,p}(\mathbb D\cap U,\mathbb C)
\cap C^0(\overline{\mathbb D}\cap U,\mathbb C)
\qquad(1\le a,j\le r).
\]
\end{rem}

\begin{proof}[Proof of (1)]

Since $p>2$, the Sobolev embedding gives
%\cite[Theorem~3.49(3)]{Ch4Hunter}より
\begin{equation}
\label{eq:ch4-sobolev-embedding}
W^{1,p}(\mathbb D)\hookrightarrow C^0(\overline{\mathbb D}),
\qquad
|u(\zeta)|\le\|u\|_\infty\le C\|u\|_{W^{1,p}}
\quad(\zeta\in\overline{\mathbb D}).
\end{equation}
Thus point evaluation is bounded.
In particular,
\[
\|uv\|_{W^{1,p}}
\underalign{\text{by $\nabla(uv)=u\nabla v+v\nabla u$}}{\le}
C\bigl(
\|u\|_\infty\|v\|_{W^{1,p}}
+\|v\|_\infty\|u\|_{W^{1,p}}
\bigr)
\underalign{\eqref{eq:ch4-sobolev-embedding}}{\le } C\|u\|_{W^{1,p}}\|v\|_{W^{1,p}}.
\]
Hence $u, v \in W^{1,p}$ implies $uv \in W^{1,p}$.

\begin{claim}
\label{claim-construction-right-inverse}
Let $p>2$ and $A\in L^p(\mathbb D;\operatorname{End}\mathbb C^m)$, and define
\[
D_Av:=\bar\partial v+Av
\]
For every $\zeta\in\overline{\mathbb D}$,
there exists a bounded complex linear map
\[
R_\zeta:L^p(\mathbb D,\mathbb C^m)
\longrightarrow
\{v\in W^{1,p}(\mathbb D,\mathbb C^m):v(\zeta)=0\}
\]
such that, for all $b\in L^p(\mathbb D,\mathbb C^m)$,
\begin{equation}
\label{eq:ch4-point-right-inverse}
D_A(R_\zeta b)=b,\qquad
(R_\zeta b)(\zeta)=0,\qquad
\|R_\zeta b\|_{W^{1,p}}\le C\|b\|_{L^p}
\end{equation}
Here $C>0$ is independent of $b$.
\end{claim}
\begin{proof}[Proof of Claim \ref{claim-construction-right-inverse}]

First suppose that $\zeta=0$.
Following \cite[\S2, immediately before Proposition~2.1]{Ch4ST},
define the Cauchy--Green operator 
$$
T:L^p(\mathbb D,\mathbb C^m)\to W^{1,p}(\mathbb D,\mathbb C^m)
\quad
(Tb)(z):=\frac1\pi\int_{\mathbb D}
\frac{b(w)}{z-w}\,dS(w)
$$
By \cite[Proposition~2.1(i)]{Ch4ST},
$T$ is bounded and linear, and satisfies $\bar\partial(Tb)=b$.

Let $w_1,\ldots,w_d$ be a real basis of the kernel of
the operator $u\mapsto u+T(Au)-(T(Au))(0)$.
Applying \cite[Theorem~3.1(i)--(ii), equation~(10)]{Ch4ST}
with $B_1=-A$ and $B_2=0$,
we obtain holomorphic vector-valued polynomials $p_j$ satisfying $p_j(0)=0$ such that the operator defined by
\[
\widetilde P:
W^{1,p}(\mathbb D,\mathbb C^m)
\longrightarrow
W^{1,p}(\mathbb D,\mathbb C^m)
\quad
\widetilde Pu
:=
u+T(Au)-(T(Au))(0)
+
\underbrace{
\sum_{j=1}^d
\left(\operatorname{Re}\int_{\mathbb D}w_j^*u\,dS\right)p_j
}_{\text{correction term}}
\]
is
bounded and invertible.
Since the correction term is holomorphic and vanishes at the origin,
\begin{equation}
\label{eq-property-widetildeu}
\bar\partial(\widetilde Pu)
=\bar\partial u+\underbrace{\bar\partial T(Au)}_{=Au
\text{ by $\bar\partial(Tb)=b$}}
=D_Au,
\qquad
(\widetilde Pu)(0)=u(0).
\end{equation}
We define $R_0:L^p(\mathbb D,\mathbb C^m)
\longrightarrow
\{v\in W^{1,p}(\mathbb D,\mathbb C^m):v(\zeta)=0\}
$
by
\[
R_0b:=
\frac12\left[
\widetilde P^{-1}\bigl(Tb-(Tb)(0)\bigr)
-\sqrt{-1}\,\widetilde P^{-1}\bigl(T(\sqrt{-1}b)-(T(\sqrt{-1}b))(0)\bigr)
\right]
\]
Since $T,\widetilde{P}$ are bounded, this is a bounded complex linear map satisfying
\[
D_AR_0b
\underalign{\eqref{eq-property-widetildeu}}{=}
\bar\partial(\widetilde PR_0b)
\underalign{\text{ by $\bar\partial(Tb)=b$}}{=}
\frac12\bigl(b-\sqrt{-1}(\sqrt{-1}b)\bigr)=b,\qquad
(R_0b)(0)=0
\]
Thus we obtain the desired right inverse $R_0$.

For a general $\zeta\in\overline{\mathbb D}$,
we have $\overline{\mathbb D}\subset B(\zeta,3)$.
Extend $A,b$ by zero to $B(\zeta,3)$ and denote the extensions by $\widehat A,\widehat b$.
The extensions still belong to $L^p$.
The construction at the origin gives a bounded complex linear right inverse $R_0$ satisfying
\[
\bigl(\partial_{\bar w}+3\widehat A(\zeta+3w)\bigr)R_0c=c,
\qquad (R_0c)(0)=0
\]
Set
\(
v(w):=R_0\!\left(3\widehat b(\zeta+3\bullet)\right)(w)
\). 
Then $v(0)=0$, and $v$ satisfies
\[
\bar\partial_w v
+3\widehat A(\zeta+3w)v
=3\widehat b(\zeta+3w)
\quad 
\text{ on $\D$.}
\]
Define
\(
(R_\zeta b)(z):=v (\frac{z-\zeta}{3})
\)
for any $z \in \D$.
The chain rule gives
\(
D_A(R_\zeta b)=b\) and \(
(R_\zeta b)(\zeta)=v(0)=0.
\). Since
zero extension, the above change of variables, and restriction are bounded and complex linear,
\(
\|R_\zeta b\|_{W^{1,p}(\mathbb D)}
\le C\|b\|_{L^p(\mathbb D)}
\)
also holds. This proves Claim \ref{claim-construction-right-inverse}.
\end{proof}

We now construct a frame by induction.

(a). The rank 1 case.
Since the closed disc is contractible, $E$ has a continuous frame $e_0$
(cf. \cite[Corollary~1.8]{Ch4Hatcher})
Approximate its local coefficients uniformly by smooth functions
and combine them using a partition of unity to obtain a $W^{1,p}$ section $e$.
A sufficiently close approximation preserves nonvanishing on the closed disc,
so $e$ is a frame. Thus,
\[
e:\overline{\mathbb D}\longrightarrow E,
\qquad
e(z)\in E_z\setminus\{0\},
\qquad
e\in W^{1,p}(\mathbb D,E).
\]

With respect to the frame $e$, write
$\bar\partial_E e=eA\,d\bar z$, where $A\in L^p(\mathbb D,\mathbb C)$,
and consider $D_A=\bar\partial+A$ acting on coefficient functions.
For any $\zeta\in\overline{\mathbb D}$ and $\xi\in E_\zeta$,
choose $c\in\mathbb C$ such that $\xi=c\,e(\zeta)$.
Using the right inverse in Claim~\ref{claim-construction-right-inverse},
define a section by
\[
u(z):=e(z)\bigl(c-R_\zeta(Ac)(z)\bigr)\in E_z
\]
The product estimate gives $u\in W^{1,p}(\mathbb D,E)$, and
\[
\bar\partial_Eu
=
e\underbrace{\bigl(Ac-D_AR_\zeta(Ac)\bigr)}
_{=0\text{ by }\eqref{eq:ch4-point-right-inverse}}\,d\bar z
=0,
\qquad
u(\zeta)=e(\zeta)c=\xi.
\]
Hence $u$ is holomorphic in $\mathbb D$ by \cite[Chapter~I, Corollary~3.30]{Ch4Demailly}.

For each $\zeta\in\overline{\mathbb D}$,
the above construction gives a global holomorphic $W^{1,p}$ section
$u_\zeta$ such that $u_\zeta(\zeta)\ne0$.
By the continuity following from \eqref{eq:ch4-sobolev-embedding},
there is a neighborhood $U_\zeta\subset\overline{\mathbb D}$ of $\zeta$
on which $u_\zeta$ is nonvanishing.
By compactness of the closed disc, choose finitely many such neighborhoods with
\[
\overline{\mathbb D}
=
U_{\zeta_1}\cup\cdots\cup U_{\zeta_N}
\]

Every point belongs to some $U_{\zeta_i}$,
and $u_{\zeta_i}\ne0$ at that point.
Thus $u_{\zeta_1},\ldots,u_{\zeta_N}$ have no common zero.
We show that, for some $c_1,\ldots,c_N\in\mathbb C$, the section
\[
s:\overline{\mathbb D}\to E,\qquad
s(z):=\sum_{i=1}^N c_i u_{\zeta_i}(z)
\quad
s\in W^{1,p}(\mathbb D,E).
\]
is nonvanishing on the entire boundary.
Write $u_{\zeta_i}=e\,a_i$.
Since there is no common zero, the boundary is covered by finitely many closed arcs
such that at least one $a_i$ is nonzero near each arc.
Fix one arc $\Gamma$ and reorder the indices so that $a_1\ne0$.
Then
\[
s(z)=0
\iff
c_1=-\sum_{i=2}^N c_i\frac{a_i(z)}{a_1(z)}
\qquad(z\in\Gamma).
\]
Fix $c_2,\ldots,c_N$.
The function
\[
z\longmapsto-\sum_{i=2}^N c_i\frac{a_i(z)}{a_1(z)}
\]
belongs to $W^{1,p}$ near $\Gamma$.
Using a cutoff and the Sobolev extension theorem
\cite[Theorem~3.41 and \S3.11]{Ch4Hunter},
we extend it, without changing its values on $\Gamma$,
to a $W^{1,p}$ function on a neighborhood on both sides of $\Gamma$.
Since $\mathcal L^2(\Gamma)=0$,
applying Lusin's $N$ property
\cite[Corollary~1, p.~791]{Ch4MarcusMizel} to this extension gives
\begin{equation}
\label{eq:ch4-boundary-image-zero}
\mathcal L^2\left(
\left\{-\sum_{i=2}^N c_i\frac{a_i(z)}{a_1(z)}
:z\in\Gamma\right\}
\right)=0.
\end{equation}
Thus the set of $c_1$ for which a zero occurs on $\Gamma$
has measure zero.

The set of all exceptional coefficients is closed by compactness of $\Gamma$,
and hence has measure zero in $\mathbb C^N$ by Fubini's theorem.
Choosing coefficients outside the exceptional sets for the finitely many arcs gives,
for every $z\in\partial\mathbb D$,
\(
s(z)\ne0
\)
.

Since $s$ is continuous and nonzero on the boundary, it is nonzero near the boundary.
Its interior zeros are isolated and therefore finite in number.
Denote them by $z_1,\ldots,z_\nu$, with orders $m_1,\ldots,m_\nu$, and set
\[
H(z):=\frac{s(z)}{\prod_{j=1}^{\nu}(z-z_j)^{m_j}}
\]
Then $H$ is holomorphic in the interior and nonzero on the closed disc,
and $H\in W^{1,p}(\mathbb D,E)$.
Thus, for every $z\in\overline{\mathbb D}$,
the map \(H(z) : \mathbb C\longrightarrow E_z \) defined by $c\longmapsto cH(z)$
is a complex linear isomorphism.
Since it is nonzero on the disc,
$H^{-1}$ is also continuous and belongs to $W^{1,p}$.

(b). General case.
Suppose that the frame has been constructed in rank $r-1$.
Repeating (a) with matrix-valued $A$, using a nonzero
component in the boundary argument,
choose a nowhere vanishing holomorphic $W^{1,p}$ section $s$.
The quotient bundle $E/\mathbb Cs$ also has local holomorphic $W^{1,p}$ frames,
so by induction it has a global holomorphic frame
$\widehat v_2,\ldots,\widehat v_{r}$.

Combine the standard Hermitian metrics in local frames using a partition of unity, and set
\[
E=\mathbb Cs\oplus(\mathbb Cs)^\perp,\qquad
\Pi(v)=v-s\frac{\langle s,v\rangle}{\langle s,s\rangle}
\]
Since $s$ is nonzero on the closed disc, the product and composition rules give $\Pi\in W^{1,p}$.
Apply $\Pi$ to local $W^{1,p}$ representatives of $\widehat v_j$.
The differences between the representatives vanish under $\Pi$, since $\Pi(s)=0$,
yielding a global $W^{1,p}$ section $v_j$.
Its image in the quotient is holomorphic, so
\[
\bar\partial_Ev_j=s\alpha_j\,d\bar z,\qquad \alpha_j\in L^p.
\]
Correcting by the right inverse $R_0$ from Claim~\ref{claim-construction-right-inverse}
with $m=1$ and $A=0$, we obtain
\[
\bar\partial_E(v_j-sR_0\alpha_j)
=
s\underbrace{\bigl(\alpha_j-\bar\partial(R_0\alpha_j)\bigr)}
_{=0\text{ by }\eqref{eq:ch4-point-right-inverse}}
\,d\bar z=0.
\]
Therefore,
\(
H=(s,v_2-sR_0\alpha_2,\ldots,v_{r}-sR_0\alpha_{r})
\)
is a holomorphic $W^{1,p}$ frame.
By the inverse matrix formula and the product and composition rules,
$H^{-1}$ is also continuous and belongs to $W^{1,p}$.
This proves (1).
\end{proof}

\begin{proof}[Proof of (2)]
Take the frame $H=(H_1,\ldots,H_r)$
from Theorem~\ref{thm:disc-frame-variation-hessian}(1).
Define the positive definite Hermitian matrix representing the metric $h$ in this frame by
\begin{equation}
\label{eq-representation-H*hH}
G(z):=\bigl(h_z(H_i(z),H_j(z))\bigr)_{i,j=1}^{r}
\qquad(z\in\overline{\mathbb D})
\end{equation}
Here the inner product is complex linear in the second variable.
We abbreviate this as $G=H^*hH$.

Since $G\in W^{1,p}(\mathbb D)$ and $p>2$,
the Sobolev embedding \cite[Theorem~3.49(3)]{Ch4Hunter} gives
\[
G\in C^{0,\alpha}(\overline{\mathbb D},M_{r}(\mathbb C)),
\qquad 0<\alpha<1-\frac2p.
\]
Apply Lempert's theorem
\cite[Th\'eor\`eme~B, p.~442; proof, pp.~473--474]{Ch4Lempert}
with boundary matrix $W$
equal to $G|_{\partial\mathbb D}$ and H\"older exponent $\alpha$.
We obtain a matrix $A$ satisfying
\begin{equation}
\label{eq:ch4-lempert-factorization}
A\in C^\alpha(\overline{\mathbb D},GL(r,\mathbb C)),
\qquad
A|_{\mathbb D}\text{ is holomorphic},
\qquad
A^*A=G\quad\text{on }\partial\mathbb D
\end{equation}

Next, using $G\in W^{1,p}(\mathbb D,M_{r}(\mathbb C))$,
we show that $A\in W^{1,p}(\mathbb D,M_{r}(\mathbb C))$.
On the disc $B(0,3/2)$ of radius $3/2$, define
\[
\widetilde A(z):=
\begin{cases}
A(z),&|z|\le1,\\[2mm]
\bigl(A(1/\bar z)^*\bigr)^{-1}G(1/\bar z),
&1<|z|<3/2
\end{cases}
\]
On the circle $|z|=1$,
\[
(A^*)^{-1}G
\underalign{\eqref{eq:ch4-lempert-factorization}}{=}
(A^*)^{-1}A^*A=A
\]
so $\widetilde A\in C^0(B(0,3/2),M_{r}(\mathbb C))$.
On the annulus $\{1<|z|<3/2\}$,
$\bigl(A(1/\bar z)^*\bigr)^{-1}$ is holomorphic and bounded, so
\begin{equation}
\label{eq:ch4-reflected-factor}
\partial_{\bar z}\widetilde A=
\begin{cases}
0,&|z|<1,\\[1mm]
\bigl(A(1/\bar z)^*\bigr)^{-1}
\partial_{\bar z}\bigl(G(1/\bar z)\bigr),
&1<|z|<3/2
\end{cases}
\in L^p(B(0,3/2),M_{r}(\mathbb C)).
\end{equation}
Here we used that the reflection $z\mapsto1/\bar z$ is
a smooth diffeomorphism from $\{1<|z|<3/2\}$ to $\{2/3<|z|<1\}$.

Taking limits in integration by parts on both sides of the circle,
the boundary terms cancel by continuity.
Thus \eqref{eq:ch4-reflected-factor} also holds
in the sense of distributions on all of $B(0,3/2)$.

Rescale the Cauchy--Green operator used in (1) and define
\[
T_{3/2}:L^p(B(0,3/2),M_{r}(\mathbb C))
\longrightarrow W^{1,p}(B(0,3/2),M_{r}(\mathbb C)),
\quad
(T_{3/2}b)(z):=\frac1\pi\int_{|w|<3/2}\frac{b(w)}{z-w}\,dS(w)
\]
Applying \cite[Proposition~2.1(i)]{Ch4ST} to each component,
we see that $T_{3/2}$ is bounded and satisfies $\partial_{\bar z}(T_{3/2}b)=b$.
Hence \eqref{eq:ch4-reflected-factor} gives
\[
T_{3/2}(\partial_{\bar z}\widetilde A)\in W^{1,p}(B(0,3/2),M_{r}(\mathbb C)),
\qquad
\partial_{\bar z}\left(\widetilde A-T_{3/2}(\partial_{\bar z}\widetilde A)\right)=0.
\]
The expression in parentheses is holomorphic on $B(0,3/2)$
by \cite[Chapter~I, Corollary~3.30]{Ch4Demailly}. Since $\overline{\mathbb D}\Subset B(0,3/2)$,
we obtain $A=\widetilde A|_{\mathbb D}\in W^{1,p}(\mathbb D,M_{r}(\mathbb C))$.
Moreover, $A^{-1}$ is bounded on the closed disc, so
\[
\partial_\beta(A^{-1})
=-A^{-1}(\partial_\beta A)A^{-1}\in L^p(\mathbb D,M_{r}(\mathbb C))
\quad(\beta=x,y),
\quad
A^{-1}\in W^{1,p}(\mathbb D,M_{r}(\mathbb C)).
\]

By the product estimate in (1), the local matrix representations of $HA^{-1}$ and its inverse $AH^{-1}$
belong to $W^{1,p}$.
The frame $HA^{-1}$ is holomorphic in the interior, continuous and invertible on the closed disc, and satisfies
\[
(HA^{-1})^*h(HA^{-1})
=(A^{-1})^*GA^{-1}
\underalign{\eqref{eq:ch4-lempert-factorization}}{=}
(A^{-1})^*\underbrace{A^*A}_{G}A^{-1}=I.
\]
Thus $HA^{-1}$ is the frame required in (2).
\end{proof}

\begin{thm}[{cf.~\cite[Proposition~4.3]{OpenAI2026Canonical}}]
\label{thm:disc-center-fixed-deformation}

Let $X$ be a smooth compact complex manifold of complex dimension $n\ge1$,
and let $2<p<\infty$.
Let $f:\mathbb D\to X$ be a holomorphic map
belonging to $W^{1,p}(\mathbb D,X)$ and continuous on the closed disc. Set
\(
x:=f(0) \in X
\).

Then for every holomorphic section $V\in W^{1,p}(\mathbb D,f^*TX)$ satisfying $V(0)=0$,
there exist $\varepsilon>0$ and a continuous map
\[
F:\{s\in\mathbb C:|s|<\varepsilon\}
\times\overline{\mathbb D}\longrightarrow X
\]
such that each $F(s,\cdot)$ is holomorphic on $\mathbb D$,
the map $s\mapsto F(s,\cdot)$ is holomorphic with values in
$W^{1,p}(\mathbb D,X)$, and
\[
F(s,0)=x,\qquad F(0,z)=f(z),\qquad
F_s(0,z)=V(z),\qquad F_{\bar s}(s,z)=0.
\]
\end{thm}

\begin{proof}
From $f\in W^{1,p}$, 
$f^*TX$ has local holomorphic $W^{1,p}$ frames continuous up to the boundary.
Apply Theorem~\ref{thm:disc-frame-variation-hessian}(1), 
we  obtain such a frame $H$ on $f^*TX$.

By \cite[Theorem~1.1(iii), \S2]{Ch4Forstneric}, the space of holomorphic $W^{1,p}$ maps
$\mathbb D\to X$ is a complex Banach manifold and its tangent space at $f$ is the space of holomorphic sections
\(
\{u\in W^{1,p}(\mathbb D,f^*TX):
\bar\partial_{f^*TX}u=0\}
\).
Define
\[
\Psi:
\{g\in W^{1,p}(\mathbb D,X):g\text{ is holomorphic}\}
\longrightarrow X,
\qquad
\Psi(g):=g(0)
\]
By \cite[end of \S2, p.~117]{Ch4Forstneric}, $\Psi$ is holomorphic,
and its differential at $f$ is
\[
d\Psi_f:
\{u\in W^{1,p}(\mathbb D,f^*TX):
\bar\partial_{f^*TX}u=0\}
\longrightarrow T_xX,
\qquad
d\Psi_f(u)=u(0).
\]
This differential has a bounded complex linear right inverse
\(
\xi\longmapsto HH(0)^{-1}\xi
\)
(Here we use that $H$ is holomorphic and belongs to $W^{1,p}$.)
Since $d\Psi_f$ has a bounded complex linear right inverse,
the complex Banach implicit function theorem gives an open neighborhood of zero
\[
U_0\subset\ker d\Psi_f
=
\{\eta\in W^{1,p}(\mathbb D,f^*TX):
\eta(0)=0,\ \bar\partial_{f^*TX}\eta=0\}
\]
and a holomorphic local parametrization of holomorphic discs with fixed center,
\[
\Phi:U_0\longrightarrow\Psi^{-1}(x),
\qquad
\Phi(0)=f,\qquad d\Phi_0(\eta)=\eta
\quad \text{$( \forall \eta \in \ker d\Psi_f)$}.
\]
Since $V(0)=0$ and $\bar\partial_{f^*TX}V=0$,
we can choose $\varepsilon'>0$ sufficiently small
so that $sV\in U_0$ for $|s|<\varepsilon'$.
Set
\[
F:\{|s|<\varepsilon'\}\longrightarrow W^{1,p}(\mathbb D,X),
\qquad
F(s,\cdot):=\Phi(sV)
\]
Each $F(s,\cdot)$ is holomorphic and satisfies
\(
F(s,0)=x\) and \(F(0,\cdot)=f\).

The holomorphicity of $\Phi$ and $d\Phi_0=\mathrm{id}$ give,
in local holomorphic coordinates on the target,
\(
F_s(0,\cdot)=V\) and \(F_{\bar s}(s,\cdot)=0.
\)
Moreover, $s\mapsto F(s,\cdot)$ is holomorphic and hence real $C^\infty$.
By \eqref{eq:ch4-sobolev-embedding},
each $F(s,\cdot)$ is continuous on the closed disc.
This proves the theorem.
\end{proof}

\section{Second Variation Formulas and an Inequality for Maximizing Discs}
\label{sec:finite-matrix-trace-unified}

\subsection{Second variation formulas}
\begin{setup}
\label{setup-second-variation}
Let $X$ be a compact K\"ahler manifold of complex dimension $n$,
let $h,h_0$ be smooth K\"ahler metrics,
and let $q$ be a smooth closed real $(1,1)$-form.
Let $2<p<\infty$,
let $f\in W^{1,p}(\mathbb D,X)$ be a holomorphic map continuous on the closed disc,
and set $x:=f(0)$.

Take the frame $B=(B_1,\ldots,B_n)$ from Theorem~\ref{thm:disc-frame-variation-hessian}(2),
which is orthonormal with respect to $h$ on the boundary.
For each $j=1,\ldots,n$,
\[
B_j:\overline{\mathbb D}\longrightarrow f^*TX,
\qquad
B_j(z)\in(f^*TX)_z\simeq T_{f(z)}X.
\]
is a holomorphic $W^{1,p}$ section.

Fix $j\in\{1,\ldots,n\}$.
By Theorem~\ref{thm:disc-center-fixed-deformation},
there exist $\varepsilon>0$ and a family of discs holomorphic in the parameter $s$,
\[
F:\{|s|<\varepsilon\}\longrightarrow W^{1,p}(\mathbb D,X),
\qquad
F(s,0)=x,\quad F(0,z)=f(z)
\]
(The dependence of $F$ on $j$ is suppressed in the notation.)
Take an open set $U\subset X$ containing $x$ and local holomorphic coordinates
$w=(w^1,\ldots,w^n):U\to w(U)\subset\mathbb C^n$.
For $z\in f^{-1}(U)$, write
\[
B_j(z)=\sum_{\alpha=1}^n B_j^\alpha(z)
\left.\frac{\partial}{\partial w^\alpha}\right|_{f(z)}
\in  T_{f(z)}X.
\]
Then
\begin{equation}
\label{eq:ch5-local-velocity}
\left.\partial_s(w^\alpha\circ F)(s,z)\right|_{s=0}
=zB_j^\alpha(z),\qquad
\partial_{\bar s}(w^\alpha\circ F)(s,z)=0
\qquad(\alpha=1,\ldots,n)
\end{equation}
Interchanging derivatives and differentiating \eqref{eq:ch5-local-velocity}, we obtain,
for each $\alpha=1,\ldots,n$,
\begin{equation}
\label{eq:ch5-center-derivatives}
\left.\partial_z\partial_s(w^\alpha\circ F)\right|_{s=0,z=0}
\underalign{\eqref{eq:ch5-local-velocity}}{=}
B_j^\alpha(0)
\quad
\left.\partial_z\partial_{\bar s}(w^\alpha\circ F)\right|_{s=0,z=0}
%\underalign{\partial_{\bar s}(w^\alpha\circ F)=0}{=}0,
=
\left.\partial_z\partial_s\partial_{\bar s}(w^\alpha\circ F)\right|_{s=0,z=0}
%\underalign{\partial_{\bar s}(w^\alpha\circ F)=0}{=}
=
0
\end{equation}

\end{setup}

\begin{thm}[{cf.~\cite[Proposition~5.2]{OpenAI2026Canonical}; cf.~\cite[Lemma~3.3, Proposition~3.4]{DuGuoXie2026}}]
\label{thm:disc-second-variation}
In Setup~\ref{setup-second-variation}, for each $j=1,\ldots,n$,
\begin{equation}
\label{eq:ch5-second-variation}
\begin{aligned}
\left.\partial_s\partial_{\bar s}d(F(s,\cdot))\right|_{s=0}
&=
|B_j(0)|_{h_0}^2,\\
\left.\partial_s\partial_{\bar s}P_q(F(s,\cdot))\right|_{s=0}
&=
\frac1{2\pi}\int_0^{2\pi}q(B_j,B_j)(e^{\sqrt{-1}\theta})\,d\theta,\\
\left.\partial_s\partial_{\bar s}A_h(F(s,\cdot))\right|_{s=0}
&=
\frac1\pi\int_{\mathbb D}\partial_z\partial_{\bar z}\bigl(|z|^2h(B_j,B_j)\bigr)\,dS.
\end{aligned}
\end{equation}

\end{thm}

\begin{proof}

The first identity follows directly:
\begin{align*}
\left.\partial_s\partial_{\bar s}d(F(s,\cdot))\right|_{s=0}
&\underalign{\text{definition of $d$}}{=}
\left.\partial_s\partial_{\bar s}
|F_z(s,0)|_{h_0,x}^2\right|_{s=0}
\underalign{\eqref{eq:ch5-center-derivatives}}{=}
|B_j(0)|_{h_0}^2.
\end{align*}

We prove the second identity in \eqref{eq:ch5-second-variation}.
In local holomorphic coordinates $w=(w^1,\ldots,w^n)$,
choose a potential $\rho_q$ with $q=\sqrt{-1}\,\partial\bar\partial\rho_q$.
Setup~\ref{setup-second-variation} and the chain rule give
\[
\left.\partial_s\partial_{\bar s}\rho_q(F(s,z))\right|_{s=0}
=
\underbrace{
\sum_{\alpha,\beta}
(\rho_q)_{\alpha\bar\beta}(f(z))
zB_j^\alpha(z)\overline{zB_j^\beta(z)}
}_{|z|^2q(B_j,B_j)(z)}
+
\underbrace{
2\operatorname{Re}\sum_\alpha
(\rho_q)_\alpha(f(z))
\left.\partial_s\partial_{\bar s}(w^\alpha\circ F)(s,z)\right|_{s=0}
}_{=0\text{: }F\text{ is holomorphic in the parameter}}.
\]
Since each disc is holomorphic,
$q(F_z,F_z)=\partial_z\partial_{\bar z}\rho_q(F)$, and hence
\begin{equation}
\label{eq:ch5-density-variation}
\begin{aligned}
\left.\partial_s\partial_{\bar s}q(F_z,F_z)\right|_{s=0}
&\underalign{\text{}}{=}
\partial_z\partial_{\bar z}
\left(
\left.\partial_s\partial_{\bar s}\rho_q(F)\right|_{s=0}
\right)
\underalign{\text{chain rule}}{=}
\partial_z\partial_{\bar z}\bigl(|z|^2q(B_j,B_j)\bigr).
\end{aligned}
\end{equation}
Note that this identity holds on the entire disc.
Therefore,
\begin{align*}
\left.\partial_s\partial_{\bar s}P_q(F(s,\cdot))\right|_{s=0}
&\underalign{\text{Definition~\ref{defn-APJ-funcional}}}{=}
\frac2\pi\int_{\mathbb D}\log\frac1{|z|}
\left.\partial_s\partial_{\bar s}q(F_z,F_z)\right|_{s=0}\,dS
\underalign{\eqref{eq:ch5-density-variation}}{=}
\frac2\pi\int_{\mathbb D}\log\frac1{|z|}
\,\partial_z\partial_{\bar z}\bigl(|z|^2q(B_j,B_j)\bigr)\,dS
\\
&\underalign{\eqref{eq-lem:green-normalization-2}}{=}\frac1{2\pi}\int_0^{2\pi}q(B_j,B_j)(e^{\sqrt{-1}\theta})\,d\theta
- \underbrace{\left.(|z|^2q(B_j,B_j))\right|_{z=0}}_{=0}
\end{align*}
Here the holomorphicity and $W^{1,p}$ regularity of $f,B_j$ give
$\partial_z\partial_{\bar z}(|z|^2q(B_j,B_j))\in L^{p/2}(\mathbb D)$,
so we  can apply \eqref{eq-lem:green-normalization-2} in Lemma~\ref{lem:green-normalization}.

For the third identity, set $q=\omega_h$ in \eqref{eq:ch5-density-variation}
and integrate over the disc with factor $1/\pi$ and without the logarithmic weight.
\end{proof}

\subsection{An inequality for locally maximizing discs}

\begin{thm}[{cf.~\cite[Corollary~5.3]{OpenAI2026Canonical}}]
\label{thm:signed-disc-determinant-trace}
Let $\delta>0$.
In Setup~\ref{setup-second-variation}, suppose that $f : \D \to X$ is a local maximum of
\(
J_\delta(f):=d(f)+P_q(f)-\delta A_h(f)
\)
among holomorphic $W^{1,p}$ discs with $f(0)=x$.
Then
\begin{equation}
\label{eq:ch5-main}
\begin{aligned}
n\exp\!\left(\frac{
P_{\operatorname{Ric}(h)}(f)-\log\frac{\omega_h^n}{\omega_{h_0}^n}(x)}n\right)
+\frac1{2\pi}\int_0^{2\pi}
(\operatorname{tr}_h q)(f(e^{\sqrt{-1}\theta}))\,d\theta
\le\delta\bigl\{n-A_{\operatorname{Ric}(h)}(f)\bigr\}.
\end{aligned}
\end{equation}

\end{thm}

\begin{proof} 
Since $f$ is a local maximum, the family $F$ in Setup~\ref{setup-second-variation}, with $F(0,z)=f(z)$, gives
\begin{equation}
\label{eq:ch5-local-maximum}
\begin{aligned}
0
&\underalign{\text{local maximality at $s=0$}}{\ge}
\left.\partial_s\partial_{\bar s}J_\delta(F(s,\cdot))\right|_{s=0}\\
&\underalign{\eqref{eq:ch5-second-variation}}{=}
|B_j(0)|_{h_0}^2
+\frac1{2\pi}\int_0^{2\pi}q(B_j,B_j)(e^{\sqrt{-1}\theta})\,d\theta
-\frac\delta\pi\int_{\mathbb D}
\partial_z\partial_{\bar z}\bigl(|z|^2h(B_j,B_j)\bigr)\,dS.
\end{aligned}
\end{equation}
For the frame $B=(B_1,\ldots,B_n)$
chosen in Setup~\ref{setup-second-variation}, write the metric $h$ as
\[
G(z):=\bigl(h_{f(z)}(B_i(z),B_j(z))\bigr)_{i,j=1}^n
=B(z)^*h_{f(z)}B(z),
\qquad
G|_{\partial\mathbb D}
\underalign{\text{orthonormality of $B$ on $\partial\mathbb D$}}{=}I
\]
(See also \eqref{eq-representation-H*hH}.)
Thus we obtain 
\begin{equation}
\label{eq:ch5-frame-identities}
\left.zB_j(z)\right|_{z=0}=0,\quad
\left.\frac d{dz}(zB_j(z))\right|_{z=0}=B_j(0),\quad
\sum_{j=1}^nh(zB_j,zB_j)=|z|^2\operatorname{tr}G,\quad
\end{equation}
and 
\begin{equation}
\label{eq:ch5-frame-identities-2}
\sum_{j=1}^nq(zB_j,zB_j)
\underalign{|z|=1}{=}
\sum_{j=1}^nq(B_j,B_j)
\underalign{\text{orthonormality of $B$ on $\partial\mathbb D$}}{=}
(\operatorname{tr}_h q)\circ f
\quad
\text{ on $z \in \partial\mathbb D$}
\end{equation}

The holomorphicity and $W^{1,p}$ regularity of $f,B$ give
$G\in W^{1,p}(\mathbb D,M_n(\mathbb C))$ and
$\Delta G\in L^{p/2}(\mathbb D,M_n(\mathbb C))$.
Summing \eqref{eq:ch5-local-maximum} over $j=1,\ldots,n$, we obtain
\begin{equation}
\label{eq:ch5-frame-sum}
\begin{aligned}
&\sum_{j=1}^n|B_j(0)|_{h_0}^2
+\frac1{2\pi}\int_0^{2\pi}
\underbrace{\sum_{j=1}^nq(B_j,B_j)(e^{\sqrt{-1}\theta})}_{= (\operatorname{tr}_h q)(f(e^{\sqrt{-1}\theta})) \text{ by \eqref{eq:ch5-frame-identities-2} }}\,d\theta
\underalign{\eqref{eq:ch5-local-maximum}}{\le} \frac\delta\pi\int_{\mathbb D}
\partial_z\partial_{\bar z}\bigl(|z|^2\sum_{j=1}^nh(B_j,B_j)\bigr)\,dS
\\
&\underalign{\eqref{eq:ch5-frame-identities}}{=} 
\frac\delta\pi\int_{\mathbb D}
\partial_z\partial_{\bar z}(|z|^2\operatorname{tr}G)\,dS
\underalign{\eqref{eq:green-radial-product}}{=}
\underbrace{\frac\delta{2\pi}\int_0^{2\pi}
\operatorname{tr}G(e^{\sqrt{-1}\theta})\,d\theta}_{=\delta n \text{ by $G|_{\partial\mathbb D} = I$}}
+\frac\delta\pi\int_{\mathbb D}
\partial_z\partial_{\bar z}\operatorname{tr}G\,dS
\end{aligned}
\end{equation}

Thus it suffices to prove the following claim.
\begin{claim}
\label{claim:ch5-curvature-and-center}
\begin{equation}
\label{eq:ch5-normalized-flux}
\frac1\pi\int_{\mathbb D}
\partial_z\partial_{\bar z}\operatorname{tr}G\,dS
=-A_{\operatorname{Ric}(h)}(f)
\end{equation}
\begin{equation}
\label{eq:ch5-center-lower-bound}
n\exp\!\left(\frac{
P_{\operatorname{Ric}(h)}(f)-\log\frac{\omega_h^n}{\omega_{h_0}^n}(x)}n\right)
\le\sum_{j=1}^n|B_j(0)|_{h_0}^2.
\end{equation}
\end{claim}
Assuming this claim,
substitute \eqref{eq:ch5-normalized-flux} and \eqref{eq:ch5-center-lower-bound} into \eqref{eq:ch5-frame-sum}
to obtain the required inequality \eqref{eq:ch5-main}.
\end{proof}

\begin{proof}[Proof of Claim \ref{claim:ch5-curvature-and-center}]
First, in local holomorphic coordinates,
$\det G=|\det B|^2\det(h_{a\bar b}\circ f)$, so 
\begin{equation}
\label{eq-ddc-logdetG}
-\sqrt{-1}\,\partial\bar\partial\log\det G
=
-\underbrace{\sqrt{-1}\,\partial\bar\partial\log|\det B|^2}_{=0\text{: }\det B\text{ is holomorphic}}
-\sqrt{-1}\,\partial\bar\partial\log\det(h_{a\bar b}\circ f)
\underalign{\text{Definition of Ric}}{=}
f^*\operatorname{Ric}(h).
\end{equation}

To prove \eqref{eq:ch5-normalized-flux}, we show that the difference between the integrals of the Laplacians
of $\operatorname{tr}G$ and $\log\det G$ is zero.
Set $u:=\operatorname{tr}G-\log\det G-n$.
The above formula for the Laplacian and $G|_{\partial\mathbb D}=I$ give
$\Delta u\in L^1(\mathbb D)$ and $u|_{\partial\mathbb D}=0$.
Taylor expansion at $I$ gives $0\le u\le C|G-I|^2$, so
\begin{equation}
\label{eq:ch5-boundary-average}
\begin{aligned}
0\le\frac1{2\pi}\int_0^{2\pi}u(re^{\sqrt{-1}\theta})\,d\theta
&\underalign{\text{Taylor expansion}}{\le}
C\int_0^{2\pi}|G(re^{\sqrt{-1}\theta})-I|^2\,d\theta
\\
&
\underalign{
G|_{\partial\mathbb D}=I\text{ and Cauchy--Schwarz}
}{\le}
C(1-r)\underbrace{
\int_{r<|z|<1}|\nabla G|^2\,dS
}_{\to0\text{: }\nabla G\in L^2(\mathbb D)}
=o(1-r).
\end{aligned}
\end{equation}
Since $\Delta u\in L^1(\mathbb D)$,
\eqref{eq:pred-3} 
gives
\begin{equation}
\label{eq:ch5-zero-flux}
\begin{aligned}
\frac1\pi\int_{\mathbb D}u_{z\bar z}\,dS
&\underalign{\eqref{eq:pred-3}}{=}
\lim_{r\uparrow1}\frac r{4\pi}\frac d{dr}
\int_0^{2\pi}u(re^{\sqrt{-1}\theta})\,d\theta
\underalign{
\text{}
}{=}
-\frac1{4\pi}\lim_{r\uparrow1}
\frac{\int_0^{2\pi}u(re^{\sqrt{-1}\theta})\,d\theta}{1-r}
\underalign{\eqref{eq:ch5-boundary-average}}{=}0.
\end{aligned}
\end{equation}
Thus the difference between the integrals of the Laplacians of $\operatorname{tr}G$ and $\log\det G$ is zero.
Using this, we obtain
\[
\frac1\pi\int_{\mathbb D}
\partial_z\partial_{\bar z}\operatorname{tr}G\,dS
\underalign{\eqref{eq:ch5-zero-flux}
}{=}
\frac1\pi\int_{\mathbb D}
\partial_z\partial_{\bar z}\log\det G\,dS
\underalign{\text{\eqref{eq-ddc-logdetG} \& Definition~\ref{defn-APJ-funcional}}
}{=}
-A_{\operatorname{Ric}(h)}(f).
\]
This is \eqref{eq:ch5-normalized-flux}.
Next, we prove \eqref{eq:ch5-center-lower-bound}.
\begin{equation}
\label{eq:pred-14}
P_{\operatorname{Ric}(h)}(f)
\underalign{ \eqref{eq-lem:green-normalization-2} \& \eqref{eq-ddc-logdetG}}{=}
\log\det G(0)
-\underbrace{
\frac1{2\pi}\int_0^{2\pi}\log\det G(e^{\sqrt{-1}\theta})\,d\theta
}_{=0\text{: }G|_{\partial\mathbb D}=I}
=\log\det G(0).
\end{equation}
Using matrix representations in local coordinates at the center,
\begin{equation}
\label{eq:ch5-center-determinant-ratio}
\frac{\det G(0)}{\det(B(0)^*h_0(x)B(0))}
=
\frac{|\det B(0)|^2\det h(x)}
{|\det B(0)|^2\det h_0(x)}
=\frac{\omega_h^n}{\omega_{h_0}^n}(x).
\end{equation}
Therefore,
\begin{align*}
\sum_{j=1}^n|B_j(0)|_{h_0}^2
&\underalign{}{=}\operatorname{tr}(B(0)^*h_0(x)B(0)) \\
&\underalign{\text{AMGM-inequality}}{\ge}
n\det(B(0)^*h_0(x)B(0))^{1/n}
\underalign{ \eqref{eq:pred-14} \& \eqref{eq:ch5-center-determinant-ratio}}{=}
n\exp\!\left(\frac{
P_{\operatorname{Ric}(h)}(f)-\log\frac{\omega_h^n}{\omega_{h_0}^n}(x)}n\right).
\end{align*}
This proves \eqref{eq:ch5-center-lower-bound}.
\end{proof}

\subsection{From Brody hyperbolicity to a scalar curvature estimate}
\begin{prop}[{cf.~\cite[Corollary~5.3]{OpenAI2026Canonical}}]
\label{prop:ricci-reward-master-inequality}
Assume that \(X\) is Brody hyperbolic, and let $h_0$ be a fixed K\"ahler metric.
Let $D_0>0$ be the constant in Lemma~\ref{lem:brody-bound}.
Let $h$ be a smooth K\"ahler metric.
Fix $x\in X$ and $\delta>0$, and let $f$ be a disc maximizing
\[
J_\delta(g):=d(g)+P_{\operatorname{Ric}(h)}(g)-\delta A_h(g)
\]
among finite-area holomorphic discs with center $x$.
Writing $R_h=\operatorname{tr}_h\operatorname{Ric}(h)$, we have
\begin{equation}
\label{eq:pred-13}
\begin{aligned}
n\exp\!\left(\frac{
P_{\operatorname{Ric}(h)}(f)-\log\frac{\omega_h^n}{\omega_{h_0}^n}(x)}n\right)
+\frac1{2\pi}\int_0^{2\pi}R_h(f(e^{\sqrt{-1}\theta}))\,d\theta
\le\delta\bigl\{n-A_{\operatorname{Ric}(h)}(f)\bigr\}
\le\delta(n+D_0).
\end{aligned}
\end{equation}
\end{prop}

\begin{proof}
By Proposition~\ref{prop:disc-maximizer}, $f$ is continuous on the closed disc
and belongs to $W^{1,p}$ for every $2<p<4$. Fix such a $p$.
For any holomorphic $W^{1,p}$ disc $g$, H\"older's inequality gives
\[
\int_{\mathbb D}|g'|_{h_0}^2\,dS
\le\pi^{1-2/p}\left(\int_{\mathbb D}|g'|_{h_0}^p\,dS\right)^{2/p}<\infty.
\]
Thus $f$ is also a local maximum of $J_\delta$ among holomorphic $W^{1,p}$ discs with fixed center.
Substitute $q=\operatorname{Ric}(h)$ into \eqref{eq:ch5-main}
of Theorem~\ref{thm:signed-disc-determinant-trace}.
Since $\operatorname{tr}_h q=R_h$, this gives the first inequality in \eqref{eq:pred-13}.
Moreover,
\[
\delta\bigl\{n-A_{\operatorname{Ric}(h)}(f)\bigr\}
\underalign{\text{Corollary~\ref{lem:signed-reward-source-bounds}}}{\le}
\delta(n+D_0)
\]
which proves the second inequality.
\end{proof}

\begin{thm}[{cf.~\cite[Proposition~7.1]{OpenAI2026Canonical}}]
\label{thm:scalar-volume-by-finite-prediction}
Assume that \(X\) is Brody hyperbolic, and let $h_0$ be a fixed K\"ahler metric.
Let $D_0>0$ be the constant in Lemma~\ref{lem:brody-bound}.
Let $h$ be a smooth K\"ahler metric on $X$ satisfying, for some $B>0$,
\begin{equation}
\label{eq-assumption-scalar}
R_h=\operatorname{tr}_h\operatorname{Ric}(h)\ge-B
\qquad\text{on }X
\end{equation}
Then, for every $x\in X$,
\begin{equation}
\label{eq:pred-16}
\frac{\omega_h^n}{\omega_{h_0}^n}(x)\ge e^{-D_0}\left(\frac nB\right)^n>0.
\end{equation}
\end{thm}

\begin{proof}
Fix $x\in X$. For each $\delta>0$,
apply Proposition~\ref{prop:disc-maximizer} with $q=\operatorname{Ric}(h)$
to choose a disc $f_\delta$ maximizing
\(
J_\delta(g):=d(g)+P_{\operatorname{Ric}(h)}(g)-\delta A_h(g)
\)
among finite-area holomorphic discs with center $x$.
\begin{equation}
\label{eq-exp-estimate}
n\exp\!\left(\frac{
P_{\operatorname{Ric}(h)}(f_\delta)-\log\frac{\omega_h^n}{\omega_{h_0}^n}(x)}n\right)
\underalign{\eqref{eq:pred-13}}{\le}
\underbrace{-
\frac1{2\pi}\int_0^{2\pi}
R_h(f_\delta(e^{\sqrt{-1}\theta}))\,d\theta
}_{\le B\text{ by }\eqref{eq-assumption-scalar}}
+
\delta(n+D_0)
\le B+\delta(n+D_0).
\end{equation}
Taking logarithms and rearranging for $\log(\omega_h^n/\omega_{h_0}^n)(x)$, we obtain
\begin{equation}
\label{eq:pred-17}
\begin{aligned}
\log\frac{\omega_h^n}{\omega_{h_0}^n}(x)
&\underalign{\eqref{eq-exp-estimate}}{\ge}
P_{\operatorname{Ric}(h)}(f_\delta)+n\log\frac n{B+\delta(n+D_0)}
\underalign{\text{Corollary~\ref{lem:signed-reward-source-bounds}}}{\ge}
-D_0+n\log\frac n{B+\delta(n+D_0)}.
\end{aligned}
\end{equation}
Since the left-hand side is independent of $\delta$, letting $\delta\downarrow0$ gives
\(
\log\frac{\omega_h^n}{\omega_{h_0}^n}(x)\ge-D_0+n\log\frac nB.
\)
Exponentiating yields \eqref{eq:pred-16}.
\end{proof}

\begin{cor}[{\cite[Proposition~6.1]{OpenAI2026Canonical}}]
\label{cor:uniform-ricci-poisson-bound}
Assume that \(X\) is Brody hyperbolic, and let $h_0$ be a fixed K\"ahler metric.
Then there is a constant $C>0$, depending only on $X,h_0$,
such that
\[
P_{\operatorname{Ric}(h_0)}(f)\le C
\]
 for every finite-area holomorphic disc $f:\mathbb D\to X$.
\end{cor}

\begin{proof}
By compactness, choose $B\ge n$
such that $R_{h_0}\ge-B$.
Let $D_0$ be the constant in Lemma~\ref{lem:brody-bound}.

Fix a finite-area holomorphic disc $f$.
For each $\delta>0$, Proposition~\ref{prop:disc-maximizer} gives
a disc $f_\delta$ maximizing
\(
J_\delta(g)
:=d(g)+P_{\operatorname{Ric}(h_0)}(g)-\delta A_{h_0}(g)
\)
among finite-area holomorphic discs with center $f(0)$.
Applying \eqref{eq:pred-13} of Proposition~\ref{prop:ricci-reward-master-inequality}
with $h=h_0$, we obtain
\[
n\exp\!\left(
\frac{P_{\operatorname{Ric}(h_0)}(f_\delta)}n
\right)
+\underbrace{\frac1{2\pi}\int_0^{2\pi}R_{h_0}(f_\delta(e^{\sqrt{-1}\theta}))\,d\theta}_{\ge-B}
\underalign{\eqref{eq:pred-13}}{\le}
\delta(n+D_0).
\]
Therefore,
\(
P_{\operatorname{Ric}(h_0)}(f_\delta)
\le n\log\frac{B+\delta(n+D_0)}{n}.
\).
Hence,
\begin{align*}
&P_{\operatorname{Ric}(h_0)}(f)
\underalign{}{=} J_\delta(f)-d(f)+\delta A_{h_0}(f)
\underalign{\text{maximality of $f_\delta$ and $d(f)\ge0$}}{\le}
J_\delta(f_\delta)+\delta A_{h_0}(f)\\
&\underalign{}{=} 
\underbrace{d(f_\delta)}_{\le D_0 \text{ by Lemma~\ref{lem:brody-bound}}}
+\underbrace{P_{\operatorname{Ric}(h_0)}(f_\delta)}_{\le n\log\frac{B+\delta(n+D_0)}{n}}
-\underbrace{\delta A_{h_0}(f_\delta)}_{\ge0}
+\delta A_{h_0}(f)
\le
D_0+n\log\frac{B+\delta(n+D_0)}n
+\delta A_{h_0}(f).
\end{align*}

Since the fixed disc $f$ has finite area,
$\delta A_{h_0}(f)\to0$.
Thus, letting $\delta\downarrow0$, we obtain
\(
P_{\operatorname{Ric}(h_0)}(f)
\le
D_0+n\log\frac Bn
\). 
The constant on the right-hand side is independent of $f$.
\end{proof}

\section{Nefness via Disc Envelopes and Bigness via the Monge--Amp\`ere Equation}
\label{sec:canonical-nef-envelope}

\begin{thm}[{cf.~\cite[Propositions~6.2, 6.4]{OpenAI2026Canonical}}]
\label{thm:canonical-nef-envelope}
Assume that \(X\) is Brody hyperbolic, and let $h_0$ be a fixed K\"ahler metric.
For $x\in X$, define
\begin{equation}
\label{eq:canonical-disc-envelope}
u(x):=
\inf\left\{
-P_{\operatorname{Ric}(h_0)}(g):
g\text{ is holomorphic from a neighborhood of }\overline{\mathbb D}\text{ to }X
\quad g(0)=x
\right\}
\end{equation}
Then $u$ is upper semicontinuous and bounded, and satisfies
\begin{equation*}
%\label{eq:canonical-envelope-positive}
-\operatorname{Ric}(h_0)+\sqrt{-1}\,\partial\bar\partial u\ge0
\end{equation*}
in the sense of currents.
In particular, $K_X$ is nef.
\end{thm}

\begin{proof}
Every disc in \eqref{eq:canonical-disc-envelope} has finite area. Thus, with the constant $C$ of Corollary~\ref{cor:uniform-ricci-poisson-bound},
\[
-C
\underalign{\text{Corollary~\ref{cor:uniform-ricci-poisson-bound}}}{\le}
u(x)
\underalign{\text{}}{\le}0.
\]
Applying Magn\'usson's disc envelope theorem \cite[Theorem~1.1]{Ch7Magnusson}
with $\omega=-2\operatorname{Ric}(h_0)$ and $\varphi=\varphi_1=\varphi_2=0$, we obtain
\begin{equation}
\begin{aligned}
\label{eq:canonical-magnusson-functional}
&\sup\{v(x):v\text{ is }\omega\text{-psh},\ v\le0\}
\\
&\underalign{\text{\cite[Theorem~1.1]{Ch7Magnusson}}}{=}
\inf
\left\{
-\mathcal{R}_{g^*\omega}(0)
\mid g\text{ is holomorphic from a neighborhood of }\overline{\mathbb D}\text{ to }X
\quad g(0)=x
\right\}.
\end{aligned}
\end{equation}
Here $\mathcal{R}_{g^*\omega}(0)$ is the Riesz potential in \cite{Ch7Magnusson}, and by definition,
\begin{equation*}
\begin{aligned}
-\mathcal{R}_{g^*(-2\operatorname{Ric}(h_0))}(0)
&\underalign{\cite[\S4, (3)]{Ch7Magnusson}}{:=}
-\frac1{2\pi}\int_{\mathbb D}\log|z|\,
\underbrace{g^*(-2\operatorname{Ric}(h_0))}
_{=-4\operatorname{Ric}(h_0)(g',g')\,dS}\\
&\underalign{}{=}\frac2\pi\int_{\mathbb D}\log|z|\,
\operatorname{Ric}(h_0)(g',g')\,dS
\underalign{\text{Definition~\ref{defn-APJ-funcional}}}{=}
-P_{\operatorname{Ric}(h_0)}(g).
\end{aligned}
\end{equation*}
Thus the disc envelope in \eqref{eq:canonical-magnusson-functional}
agrees with $u$ in \eqref{eq:canonical-disc-envelope}.

By \cite[Theorem~1.1]{Ch7Magnusson}, $u\not\equiv-\infty$, $u$ is upper semicontinuous, and
$-\operatorname{Ric}(h_0)+\sqrt{-1}\,\partial\bar\partial u\ge0$
in the sense of currents.
In particular, $K_X$ admits a semipositive singular metric with zero Lelong numbers everywhere, and is therefore nef by \cite[Proposition~18.8(a)]{DemaillyAnalyticMethods}.
\end{proof}

\begin{thm}[{cf.~\cite[Proposition~8.2, Lemma~9.1]{OpenAI2026Canonical}, 
\cite[Proofs of Lemmas~5 and~6]{Ch8WuYau}
 }]
\label{thm:canonical-big-volume}
Assume that \(X\) is Brody hyperbolic, and let $h_0$ be a fixed K\"ahler metric.
Let $D_0$ be the constant in Lemma~\ref{lem:brody-bound}.
Then
\begin{equation}
(2\pi)^n c_1(K_X)^n
\ge e^{-D_0}\int_X \omega_{h_0}^n>0.
\end{equation}
In particular, $K_X$ is ample.
\end{thm}

\begin{proof}
By Theorem~\ref{thm:canonical-nef-envelope}, $K_X$ is nef.
Applying Wu--Yau's construction of metrics \cite[Proposition~8(i)]{Ch8WuYau},
for every $\varepsilon>0$ we obtain a smooth real function $u_\varepsilon$ such that
\begin{equation}
\omega_{h_\varepsilon}:=\varepsilon\omega_{h_0}-\operatorname{Ric}(h_0)+\sqrt{-1}\,\partial\bar\partial u_\varepsilon>0,
\quad \text{and} \quad
\omega_{h_\varepsilon}^n=e^{u_\varepsilon}\omega_{h_0}^n
\end{equation}
Then
$\operatorname{Ric}(h_\varepsilon)=\operatorname{Ric}(h_0)-\sqrt{-1}\,\partial\bar\partial u_\varepsilon
=
-\omega_{h_\varepsilon}+\varepsilon\omega_{h_0}
$, so
\begin{align*}
R_{h_\varepsilon}
=\operatorname{tr}_{h_\varepsilon}(-\omega_{h_\varepsilon}+\varepsilon\omega_{h_0})
=-n+\underbrace{\varepsilon\operatorname{tr}_{h_\varepsilon}h_0}_{\ge0}
\ge-n.
\end{align*}
Applying Theorem~\ref{thm:scalar-volume-by-finite-prediction}
with $h=h_\varepsilon$ and $B=n$, we obtain
$
\omega_{h_\varepsilon}^n
\ge
e^{-D_0}\omega_{h_0}^n
$.
On the other hand, $[\omega_{h_\varepsilon}]=2\pi c_1(K_X)+\varepsilon[\omega_{h_0}]$, so
\begin{align*}
(2\pi)^n c_1(K_X)^n
&=\lim_{\varepsilon\downarrow0}\int_X
\bigl(-\operatorname{Ric}(h_0)+\varepsilon\omega_{h_0}\bigr)^n
=\lim_{\varepsilon\downarrow0}\int_X \omega_{h_\varepsilon}^n
\underalign{\text{(by $\omega_{h_\varepsilon}^n
\ge
e^{-D_0}\omega_{h_0}^n$)}}{\ge}
e^{-D_0}\int_X \omega_{h_0}^n>0.
\end{align*}
Hence $c_1(K_X)^n>0$. Since $K_X$ is nef,
\cite[Corollary~6.19, p.~48]{DemaillyAnalyticMethods}
implies that $K_X$ is big.

Brody hyperbolicity implies that $X$ contains no rational curves,
and the bigness of $K_X$ and the K\"ahler property imply that $X$ is projective.
The argument for a smooth projective variety whose canonical bundle is nef and big
and which contains no rational curves
in \cite[p.~219]{DebarreHigherDimensional} therefore shows that
$K_X$ is ample
(cf.~\cite[Proof of Lemma~5]{Ch8WuYau}).
\end{proof}

\section{Optimistic Outlook}
\label{sec:outlook-second-variation}

\subsection{Constructing rational curves and dense entire curves}

In the construction of Du--Guo--Xie
\cite[\S1]{DuGuoXie2026},
deformations that decrease the area of holomorphic discs on a Fano manifold
are combined with enlargement of the domain to produce a rational curve.
Surprisingly, an argument similar to that of \cite{DuGuoXie2026} also appears in this paper.

\begin{question}
\label{que:outlook-disc-construction}
Can the second variation method of this paper or \cite{DuGuoXie2026} be used to construct dense entire curves? 
\end{question}
Here an entire curve $f : \C \to X$ is called dense if its image is dense in $X$ in the usual topology.

For a smooth projective variety with non-nef $K_X$,
an algebraic proof of the existence of rational curves using jets
is given in \cite[Theorem~1.1]{DongGaoGuoSunWuXie2026}.
Regarding an attempt to extend the second variation method
from the Fano case to the non-nef case, that paper reports that the exploration reached an impasse and was abandoned
(\cite[Appendix~B, p.~28 and Table~4 (R1)]{DongGaoGuoSunWuXie2026})

Constructing rational curves and entire curves is difficult.
We ask whether analytic methods can make their construction easier.
This leads to the following question.

\begin{question}
Can one develop a general framework for constructing rational curves and entire curves
analytically, using second variations and the methods in
\cite{DuGuoXie2026} and
\cite[Theorem~1.1]{DongGaoGuoSunWuXie2026}?
\end{question}

\subsection{The Oka-1 property and Campana specialness}

The Oka-1 property requires that maps from open Riemann surfaces
can be deformed to holomorphic maps within a prescribed homotopy class,
with Runge approximation and interpolation of finite jets on discrete sets
\cite[Definition~1.1]{OutlookAlarconForstneric}.
In particular, Oka-1 manifolds admit dense entire curves.

The following conjecture concerns Campana specialness
and entire curves.
\begin{conjecture}[{\cite[Conjecture~3.2]{OutlookCampanaWinkelmann}}]
\label{conj:outlook-special-entire}
A variety $X$ is special in the sense of Campana if and only if
there exists an entire curve $\mathbb C\to X$ whose image is dense in $X$.
\end{conjecture}

There is also a conjecture relating the Oka-1 property to specialness.
\begin{conjecture}[{\cite[Conjecture~10.1]{OutlookAnDuGuoWangXie}}]
\label{conj:outlook-special-oka}
For a smooth complex projective variety $X$,
does specialness imply the Oka-1 property?
\end{conjecture}
One may ask whether this can be proved using the second variation method of this paper, \cite{DuGuoXie2026}, or \cite{OutlookAnDuGuoWangXie}.

\subsection{Second variations in several variables and measure hyperbolicity}

On the hyperbolic side, there is the following conjecture.

\begin{conjecture}[One direction of the general type conjecture]
\label{conj:outlook-measure-big}
If $X$ is projective and measure hyperbolic, then $K_X$ is big
\cite[Conjecture~4.6]{OutlookDemaillyGeneralType}.
Here measure hyperbolic means that the Kobayashi--Eisenman pseudovolume form
is nondegenerate outside a proper algebraic subset.
\end{conjecture}
Can the second variation method for holomorphic discs be extended
to holomorphic maps in several variables and estimates for their Jacobians?
This suggests the possibility of deriving the bigness of the canonical bundle from measure hyperbolicity.

\subsection{Oka vs Kobayashi}

Let $n:=\dim X$.
A manifold $X$ is strongly dominable if, for every $x\in X$,
there exists a holomorphic map
\[
F_x:\mathbb C^n\longrightarrow X,
\qquad F_x(0)=x,\qquad
dF_x(0):\mathbb C^n\longrightarrow T_xX
\text{ is an isomorphism}.
\]
The preprint of An--Du--Guo--Wang--Xie proves
the equivalence between the Oka property and strong dominability
\cite[Theorem~1.4]{OutlookAnDuGuoWangXie}.
Theorem~1.3 of the same paper characterizes the Oka property
by an approximation condition for local families of holomorphic discs.

The following table organizes these questions by curves versus volume, and the Oka side versus the hyperbolic side.
It compares directions of research;
it does not assert that every manifold belongs to one of the two sides.

\begin{center}
\renewcommand{\arraystretch}{1.5}
\begin{tabular}{|p{0.15\linewidth}|p{0.45\linewidth}|p{0.45\linewidth}|}
\hline
& Oka side (positive curvature)
& Kobayashi side (negative curvature)\\
\hline
Curves and discs
&
Construction of rational curves on Fano manifolds
\cite{DuGuoXie2026}
Dense entire curves, the Oka-1 property, etc. Maps from $\C$.
&
Ampleness of $K_X$ from Kobayashi hyperbolicity.
Maps from discs.
\\
\hline
$\C^n$ and balls
&
The Oka property and strong dominability
\cite[Theorem~1.4]{OutlookAnDuGuoWangXie}.
Maps from $\C^n$.
&
Bigness of $K_X$ from measure hyperbolicity?
(Conjecture~\ref{conj:outlook-measure-big})
Maps from balls
\\
\hline
\end{tabular}
\end{center}
As in Question~\ref{que:outlook-disc-construction}, it is interesting that the arguments for positive and negative curvature are similar.
This leads to the following admittedly vague and optimistic question.
\begin{question}
Is there a theory unifying Oka theory and Kobayashi theory?
\end{question}

\appendix
\section{A Proof Using the K\"ahler--Ricci Flow}
\label{sec:krf-appendix}
Following Nomura's approach in \cite[Theorems~1.1--1.2]{Nomura}, this section proves the results of Section~\ref{sec:canonical-nef-envelope} using the K\"ahler--Ricci flow. This does not, however, simplify the argument.

\begin{thm}[Nefness and bigness of the canonical bundle via K\"ahler--Ricci Flow]
\label{thm:krf-nef-big}
Assume that \(X\) is Brody hyperbolic, and let $h_0$ be a fixed K\"ahler metric.
Let $D_0>0$ be the constant in Lemma~\ref{lem:brody-bound}.
Then the following hold.
\begin{enumerate}
\item The normalized K\"ahler--Ricci flow starting from $h_0$
exists for all $t\ge0$, and $K_X$ is nef.
\item $c_1(K_X)^n>0$.
\end{enumerate}
In particular, $K_X$ is ample.
\end{thm}
\begin{proof}
We first prepare the K\"ahler--Ricci flow argument used in both (1) and (2).
Consider the unnormalized K\"ahler--Ricci flow
\cite[(3.55), Theorem~3.2.1, p.~107]{KRF-SongWeinkove}
on its maximal interval of existence $[0,S)$:
\[
\partial_s\omega_{\widetilde h_s}
=-\operatorname{Ric}(\widetilde h_s),
\qquad
\omega_{\widetilde h_s}|_{s=0}=\omega_{h_0}.
\]
By \cite[Theorem~3.3.1, (3.137), p.~122]{KRF-SongWeinkove},
\[
S=\sup\left\{
s>0:[\omega_{h_0}]+2\pi s\,c_1(K_X)\text{ is a K\"ahler class}
\right\}.
\]

Following \cite[p.~128, immediately after (3.166)]{KRF-SongWeinkove},
set $T:=\log(1+S)$ and
$h_t:=e^{-t} \widetilde h_{e^t-1}$ for $0\le t<T$.
Then $h_t$ satisfies the normalized K\"ahler--Ricci flow
on its maximal interval $[0,T)$:
\[
\partial_t\omega_{h_t}=-\operatorname{Ric}(h_t)-\omega_{h_t},
\qquad\omega_{h_t}|_{t=0}=\omega_{h_0}.
\]
Taking cohomology classes in the flow equation gives
\begin{equation}
\label{eq:krf-class}
[\omega_{h_t}]
\underalign{\text{standard ODE}}{=}
e^{-t}[\omega_{h_0}]+(1-e^{-t})\,2\pi c_1(K_X).
\end{equation}
The characterization of $S$ and the identity $T=\log(1+S)$ give
\begin{equation}
\label{eq:krf-class-time}
\begin{aligned}
T
=
%\sup\left\{
%t>0:[\omega_{h_0}]+2\pi(e^t-1)c_1(K_X)
%\text{ is a K\"ahler class}
%\right\}\\
%&\underalign{\text{multiplication by }e^{-t}>0\text{ preserves the K\"ahler property}}{=}
\sup\left\{
t>0:e^{-t}[\omega_{h_0}]+(1-e^{-t})\,2\pi c_1(K_X)
\text{ is a K\"ahler class}
\right\}.
\end{aligned}
\end{equation}

The scalar curvature evolution equation
\cite[Theorem~3.2.2, (3.56), $\nu=1$]{KRF-SongWeinkove}
gives
\begin{align*}
(\partial_t-\Delta_{h_t})\bigl(e^t(R_{h_t}+n)\bigr)
&\underalign{\text{\cite[Theorem~3.2.2, (3.56)]{KRF-SongWeinkove}}}{=}
e^t\bigl(|\operatorname{Ric}(h_t)|_{h_t}^2+2R_{h_t}+n\bigr)
\underalign{\text{Cauchy--Schwarz}}{\ge}
e^t\underbrace{\left(\frac{R_{h_t}^2}{n}+2R_{h_t}+n\right)}_{
=\frac{1}{n}(R_{h_t}+n)^2
}
\ge0.
\end{align*}
Here we use $\Delta_{h_t}( \bullet):=\operatorname{tr}_{h_t}(\sqrt{-1}\,\partial\bar\partial \bullet)$.
Since $X$ is compact without boundary, the maximum principle
\cite[Proposition~3.1.7]{KRF-SongWeinkove} gives
 \(
e^t(R_{h_t}+n)
\ge
\inf_X(R_{h_0}+n)
\).
Therefore, setting $C:=\max\{0,-\inf_XR_{h_0}-n\}$, we obtain
\begin{equation}
\label{eq:krf-scalar-lower}
R_{h_t}\ge-n-Ce^{-t}.
\end{equation}
By \cite[(3.65), $\nu=1$]{KRF-SongWeinkove} and
\eqref{eq:krf-scalar-lower},
\[
\log\frac{\omega_{h_t}^n}{\omega_{h_0}^n}
\quad
\underalign{\text{\cite[(3.65)]{KRF-SongWeinkove}}}{=}
\quad
\int_0^t(-R_{h_s}-n)\,ds
\underalign{\eqref{eq:krf-scalar-lower}}{\le}
C\underbrace{\int_0^t e^{-s}\,ds}_{\le1}
\le C.
\]
On the other hand,
\begin{equation*}
\frac{\omega_{h_t}^n}{\omega_{h_0}^n}
\quad
\underalign{
\substack{
\text{Theorem~\ref{thm:scalar-volume-by-finite-prediction}}\\
B=n+Ce^{-t}
}
}{\ge}
\quad
e^{-D_0}\left(\frac{n}{n+Ce^{-t}}\right)^n
\ge e^{-D_0}\left(\frac{n}{n+C}\right)^n.
\end{equation*}
Hence,
\begin{equation}
\label{eq:krf-short-density}
\left|\log\frac{\omega_{h_t}^n}{\omega_{h_0}^n}\right|
\underalign{}{\le}
\underbrace{
\max\left\{C,\ D_0+n\log\frac{n+C}{n}\right\}
}_{=:M\text{: independent of time}}.
\end{equation}

(1). 
Suppose that $T<\infty$ and seek a contradiction.
Let $A\ge0$ be an upper bound for the holomorphic bisectional curvature of $h_0$, and set
\begin{equation}
\label{eq:krf-auxiliary-definition}
Q:=(e^t-1)\log\frac{\omega_{h_t}^n}{\omega_{h_0}^n}
-\int_0^t e^s\log\frac{\omega_{h_s}^n}{\omega_{h_0}^n}\,ds-nt
\end{equation}

\begin{claim}
\label{claim:krf-trace-auxiliary}
We have $Q(\cdot,0)=0$, $\inf_{X\times[0,T)}Q>-\infty$, and
\begin{equation}
\label{eq:krf-short-trace}
\log\operatorname{tr}_{h_t}h_0+AQ-t
\le
\log n.
\end{equation}
%\[
%(\partial_t-\Delta_{h_t})
%\bigl(\log\operatorname{tr}_{h_t}h_0+AQ-t\bigr)\le0
%\]
%を満たす。
\end{claim}

\begin{proof}[Proof of Claim \ref{claim:krf-trace-auxiliary}]
By the definition of $Q$, we obtain $Q(\cdot,0)=0$ and
\[
Q
\underalign{\eqref{eq:krf-auxiliary-definition} \& \eqref{eq:krf-short-density}}{\ge}
-M(e^t-1)
-M\underbrace{\int_0^t e^s\,ds}_{=e^t-1}
-nt
\ge-2M(e^T-1)-nT>-\infty.
\]
So it is enough to show that \eqref{eq:krf-short-trace}.
The flow equation gives
\[
\partial_t(e^t\omega_{h_t})
\underalign{\text{by flow equation}}{=}
-e^t\operatorname{Ric}(h_t)
\underalign{\text{}}{=}
-e^t\operatorname{Ric}(h_0)
+e^t\sqrt{-1}\,\partial\bar\partial\log\frac{\omega_{h_t}^n}{\omega_{h_0}^n}.
\]
Integrating in time and rearranging, we obtain
\begin{equation}
\label{eq:krf-potential}
\sqrt{-1}\,\partial\bar\partial\int_0^t e^s\log\frac{\omega_{h_s}^n}{\omega_{h_0}^n}\,ds
=e^t\omega_{h_t}-\omega_{h_0}+(e^t-1)\operatorname{Ric}(h_0).
\end{equation}
Therefore, we obtain
\begin{align*}
\partial_tQ
&\underalign{\eqref{eq:krf-auxiliary-definition}}{=}
(e^t-1)\partial_t\log\frac{\omega_{h_t}^n}{\omega_{h_0}^n}-n
%\quad\underalign{\text{(by $\partial_t\log(\omega_{h_t}^n/\omega_{h_0}^n)=-R_{h_t}-n$)}}{=}\quad
\underalign{}{=}
-(e^t-1)R_{h_t}-ne^t,\\
\sqrt{-1}\,\partial\bar\partial Q
&\underalign{
\eqref{eq:krf-auxiliary-definition}
}{=}
(e^t-1)\bigl(\operatorname{Ric}(h_0)-\operatorname{Ric}(h_t)\bigr)
-\sqrt{-1}\,\partial\bar\partial\int_0^t e^s\log\frac{\omega_{h_s}^n}{\omega_{h_0}^n}\,ds
\\
&\underalign{\eqref{eq:krf-potential}}{=}
\omega_{h_0}-e^t\omega_{h_t}-(e^t-1)\operatorname{Ric}(h_t).
\end{align*}
Hence, 
\begin{equation}
\label{eq:krf-Q-evolution}
\begin{aligned}
(\partial_t-\Delta_{h_t})Q
&\underalign{\text{ by $\Delta_{h_t}Q=\operatorname{tr}_{h_t}(\sqrt{-1}\,\partial\bar\partial Q)$}}{=}
-(e^t-1)R_{h_t}-ne^t
-
\bigl(\operatorname{tr}_{h_t}h_0-ne^t-(e^t-1)R_{h_t}\bigr)
=-\operatorname{tr}_{h_t}h_0.
\end{aligned}
\end{equation}

We use Song--Tian's parabolic Schwarz lemma \cite{ST07}
(the parabolic version of Yau's Schwarz lemma \cite{Yau78})
in the form of \cite[Proposition~2.6]{Nomura}.
Substituting $\omega_{h_t}$ for $\omega_t$ and $\omega_{h_0}$ for $\widehat\omega$,
and writing $h_t=(g_{i\bar j}(t))$ in local coordinates, we obtain
\begin{equation}
\label{eq:krf-Q-evolution-2}
(\partial_t-\Delta_{h_t})\log\operatorname{tr}_{h_t}h_0
\underalign{\text{\cite[Proposition~2.6]{Nomura}}}{\le}
1+\frac{
\sum_{i,j,k,l}g^{i\bar j}(t)g^{k\bar l}(t)
R^{h_0}_{i\bar j k\bar l}
}{\operatorname{tr}_{h_t}h_0}
\underalign{\text{by bisectional curvature}\le A}{\le}
1+\frac{A(\operatorname{tr}_{h_t}h_0)^2}
{\operatorname{tr}_{h_t}h_0}
=1+A\operatorname{tr}_{h_t}h_0.
\end{equation}
Thus,
\[
(\partial_t-\Delta_{h_t})
\bigl(\log\operatorname{tr}_{h_t}h_0+AQ-t\bigr)
\underalign{\text{\eqref{eq:krf-Q-evolution}\& \eqref{eq:krf-Q-evolution-2}}}{\le}
\underbrace{
1+A\operatorname{tr}_{h_t}h_0
-A\operatorname{tr}_{h_t}h_0-1
}_{=0}.
\]
The maximum principle \cite[Proposition~3.1.7]{KRF-SongWeinkove} gives
\begin{equation*}
\log\operatorname{tr}_{h_t}h_0+AQ-t
\underalign{\cite[Proposition~3.1.7]{KRF-SongWeinkove}}{\le}
\underbrace{\log\operatorname{tr}_{h_0}h_0+AQ(0)}_{=\log n}.
\end{equation*}
\end{proof}
We return to the proof of (1).
By Claim~\ref{claim:krf-trace-auxiliary},
$Q\ge-K_T$ for some $K_T\ge0$.
Therefore,
\begin{equation*}
\log\operatorname{tr}_{h_t}h_0
\underalign{Claim~\ref{claim:krf-trace-auxiliary} \eqref{eq:krf-short-trace}}{\le}
\log n+t-AQ
\underalign{\text{ (by $Q\ge-K_T$ \& $A\ge0$ \& $ t<T$)}}{\le}
\log n+T+AK_T.
\end{equation*}
Hence,
\begin{equation}
\label{eq-trace-h0-eval}
h_0
\underalign{\text{}}{\le}
(\operatorname{tr}_{h_t}h_0)h_t
\le ne^{T+AK_T}h_t.
\end{equation}
Since $X$ is compact, there exists $C>0$ such that
\(
-C\omega_{h_0}
\le \omega_{h_0}+\operatorname{Ric}(h_0)
\le C\omega_{h_0}.
\)
Taking $t<T$ sufficiently close to $T$ such that
$T-t\le e^{-T-AK_T}/(2nC)$, we obtain
\[
\begin{aligned}
\omega_{h_t}
+(e^{-T}-e^{-t})\bigl(\omega_{h_0}+\operatorname{Ric}(h_0)\bigr)
&\underalign{\eqref{eq-trace-h0-eval}}{\ge}
\left(
\frac{e^{-T-AK_T}}{n}
-C|e^{-T}-e^{-t}|
\right)\omega_{h_0}
\underalign{}{\ge} 
\frac{e^{-T-AK_T}}{2n}\omega_{h_0}>0.
\end{aligned}
\]
Fix such a $t<T$. Since the above form has a uniform positive lower bound,
for sufficiently small $\varepsilon>0$ we still have
\[
0<
\omega_{h_t}
+(e^{-(T+\varepsilon)}-e^{-t})
\bigl(\omega_{h_0}+\operatorname{Ric}(h_0)\bigr)
 \underalign{\eqref{eq:krf-class}}{\in} e^{-(T+\varepsilon)}[\omega_{h_0}]
+(1-e^{-(T+\varepsilon)})\,2\pi c_1(K_X)
\]
Thus $T+\varepsilon$ belongs to the set in \eqref{eq:krf-class-time},
contradicting the fact that $T$ is its supremum.
Hence $T=\infty$, and finally $K_X$ is nef by \eqref{eq:krf-class} again.

As for (2), we have
\begin{align*}
(2\pi)^n c_1(K_X)^n
&\underalign{
\eqref{eq:krf-class}
}{=}
\lim_{t\to\infty}\int_X\omega_{h_t}^n
\underalign{
\text{Theorem~\ref{thm:scalar-volume-by-finite-prediction}
and \eqref{eq:krf-scalar-lower}}
}{\ge}
e^{-D_0}
\underbrace{
\lim_{t\to\infty}
\left(\frac{n}{n+Ce^{-t}}\right)^n
}_{=1}
\int_X\omega_{h_0}^n
=e^{-D_0}\int_X\omega_{h_0}^n>0.
\end{align*}

The bigness and then the ampleness of $K_X$ follow by the same argument as at the end of Theorem~\ref{thm:canonical-big-volume}.
\end{proof}

\bibliographystyle{alpha}
\bibliography{ref_Kobayashi_KRF_20261126_en}

\newcommand{\etalchar}[1]{$^{#1}$}
\begin{thebibliography}{DGG{\etalchar{+}}26}

\bibitem[ADG{\etalchar{+}}26]{OutlookAnDuGuoWangXie}
S.~An, Yun-Heng Du, Bin Guo, P.-C. Wang, and Song-Yan Xie.
\newblock {A Local Approximation Characterization of Oka Manifolds}, 2026.
\newblock Preprint, arXiv:2610.10475v3.
  \url{https://arxiv.org/abs/2610.10475v3}.

\bibitem[AF25]{OutlookAlarconForstneric}
Antonio Alarc\'on and Franc Forstneri\v{c}.
\newblock {Oka-1 manifolds}.
\newblock {\em Math. Z.}, 311, 2025.
\newblock Paper No.~33.

\bibitem[CW23]{OutlookCampanaWinkelmann}
Fr\'ed\'eric Campana and J\"org Winkelmann.
\newblock {Dense entire curves in rationally connected manifolds}.
\newblock {\em Algebr. Geom.}, 10(5):521--553, 2023.
\newblock With an appendix by J\'anos Koll\'ar.

\bibitem[Deb01]{DebarreHigherDimensional}
Olivier Debarre.
\newblock {\em {Higher-Dimensional Algebraic Geometry}}.
\newblock Universitext. Springer-Verlag, New York, 2001.

\bibitem[Dem12a]{DemaillyAnalyticMethods}
Jean-Pierre Demailly.
\newblock {\em {Analytic Methods in Algebraic Geometry}}, volume~1 of {\em
  Surveys of Modern Mathematics}.
\newblock International Press and Higher Education Press, 2012.

\bibitem[Dem12b]{Ch4Demailly}
Jean-Pierre Demailly.
\newblock {Complex Analytic and Differential Geometry}, 2012.
\newblock Version of June 21, 2012.
  \url{https://www-fourier.univ-grenoble-alpes.fr/~demailly/manuscripts/agbook.pdf}.

\bibitem[Dem12c]{Demailly2012Hyperbolic}
Jean-Pierre Demailly.
\newblock {Hyperbolic algebraic varieties and holomorphic differential
  equations}.
\newblock {\em Acta Math. Vietnam.}, 37(4):441--512, 2012.

\bibitem[Dem20]{OutlookDemaillyGeneralType}
Jean-Pierre Demailly.
\newblock {Recent results on the Kobayashi and Green--Griffiths--Lang
  conjectures}.
\newblock {\em Jpn. J. Math.}, 15(1):1--120, 2020.

\bibitem[DGG{\etalchar{+}}26]{DongGaoGuoSunWuXie2026}
Bin Dong, Guoxiong Gao, Bin Guo, Zeming Sun, Bin Wu, and Song-Yan Xie.
\newblock {Constructing Rational Curves via Jets on Projective Varieties with
  Non-Nef Canonical Bundle}, 2026.
\newblock arXiv:2609.28465.

\bibitem[DGX26]{DuGuoXie2026}
Yun-Heng Du, Bin Guo, and Song-Yan Xie.
\newblock {Analytic Construction of Rational Curves on Fano Manifolds}, 2026.
\newblock arXiv:2609.11612v2. \url{https://arxiv.org/abs/2609.11612v2}.

\bibitem[For07]{Ch4Forstneric}
Franc Forstneri\v{c}.
\newblock {Manifolds of holomorphic mappings from strongly pseudoconvex
  domains}.
\newblock {\em Asian J. Math.}, 11(1):113--126, 2007.

\bibitem[Hat17]{Ch4Hatcher}
Allen Hatcher.
\newblock {Vector Bundles and K-Theory}, 2017.
\newblock Version 2.2, November 2017.
  \url{https://pi.math.cornell.edu/~hatcher/VBKT/VB.pdf}.

\bibitem[Hun]{Ch4Hunter}
John~K. Hunter.
\newblock {Lecture Notes on PDEs}.
\newblock Chapter~3: Sobolev spaces.
  \url{https://www.math.ucdavis.edu/~hunter/pdes/ch3.pdf}.

\bibitem[Lem81]{Ch4Lempert}
L\'aszl\'o Lempert.
\newblock {La m\'etrique de Kobayashi et la repr\'esentation des domaines sur
  la boule}.
\newblock {\em Bull. Soc. Math. France}, 109:427--474, 1981.

\bibitem[Mag12]{Ch7Magnusson}
Benedikt~Steinar Magn\'usson.
\newblock {Extremal $\omega$-plurisubharmonic functions as envelopes of disc
  functionals: generalization and applications to the local theory}.
\newblock {\em Math. Scand.}, 111(2):296--319, 2012.

\bibitem[MM73]{Ch4MarcusMizel}
M.~Marcus and V.~J. Mizel.
\newblock {Transformations by functions in Sobolev spaces and lower
  semicontinuity for parametric variational problems}.
\newblock {\em Bull. Amer. Math. Soc.}, 79(4):790--795, 1973.

\bibitem[Nom18]{Nomura}
Ryosuke Nomura.
\newblock {K\"ahler manifolds with negative holomorphic sectional curvature,
  K\"ahler--Ricci flow approach}.
\newblock {\em Int. Math. Res. Not. IMRN}, 2018(21):6611--6616, 2018.
\newblock Numbering follows arXiv:1610.01976v2.

\bibitem[{Ope}]{OpenAIMathCatalogue}
{OpenAI}.
\newblock {Mathematics manuscript collection}.
\newblock Item~051: Kobayashi's canonical-ampleness conjecture.
  \url{https://github.com/openai/math/blob/main/CONTENTS.md} (accessed November
  7, 2026).

\bibitem[{Ope}26]{OpenAI2026Canonical}
{OpenAI}.
\newblock {Canonical ampleness of compact hyperbolic K\"ahler manifolds}, 2026.
\newblock Preprint, September 23, 2026.
  \url{https://github.com/openai/math/blob/main/preprints/Canonical-ampleness-of-compact-hyperbolic-Kahler-manifolds-September-23-2026/canonical-ampleness.pdf}.

\bibitem[ST07]{ST07}
Jian Song and Gang Tian.
\newblock {The K\"ahler--Ricci flow on surfaces of positive Kodaira dimension}.
\newblock {\em Invent. Math.}, 170(3):609--653, 2007.

\bibitem[ST12]{Ch4ST}
Alexandre Sukhov and Alexander Tumanov.
\newblock {Deformations and transversality of pseudo-holomorphic discs}.
\newblock {\em J. Anal. Math.}, 116(1):1--16, 2012.
\newblock Numbering follows arXiv:0909.4899v3.

\bibitem[SW13]{KRF-SongWeinkove}
Jian Song and Ben Weinkove.
\newblock {An Introduction to the K\"ahler--Ricci Flow}.
\newblock In S.~Boucksom, P.~Eyssidieux, and V.~Guedj, editors, {\em {An
  Introduction to the K\"ahler--Ricci Flow}}, volume 2086 of {\em Lecture Notes
  in Mathematics}, pages 89--188. Springer, 2013.

\bibitem[TY17]{IntroTosattiYang}
Valentino Tosatti and Xiaokui Yang.
\newblock {An extension of a theorem of Wu--Yau}.
\newblock {\em J. Differential Geom.}, 107(3):573--579, 2017.

\bibitem[WY16]{Ch8WuYau}
Damin Wu and Shing-Tung Yau.
\newblock {Negative holomorphic curvature and positive canonical bundle}.
\newblock {\em Invent. Math.}, 204(2):595--604, 2016.
\newblock Numbering follows the arXiv version, arXiv:1505.05802.

\bibitem[Yau78]{Yau78}
Shing-Tung Yau.
\newblock {A general Schwarz lemma for K\"ahler manifolds}.
\newblock {\em Amer. J. Math.}, 100(1):197--203, 1978.

\bibitem[Zha26]{Zhang2026VolumeSchwarz}
Shiyu Zhang.
\newblock {Remark on Kobayashi's conjecture: a volume Schwarz lemma}, 2026.
\newblock Preprint, arXiv:2611.26524. \url{https://arxiv.org/abs/2611.26524}.

\end{thebibliography}
\end{document}